\documentclass[a4paper,11pt,reqno]{amsart}

\usepackage{latexsym,dsfont,amssymb,amsmath,amsthm,amscd}

\usepackage[all]{xy}

\chardef\bslash=`\\ 

\numberwithin{equation}{section}

\newtheorem{theorem}{Theorem}[section]
\newtheorem{corollary}[theorem]{Corollary}
\newtheorem{lemma}[theorem]{Lemma}
\newtheorem{proposition}[theorem]{Proposition}
\theoremstyle{remark}
\newtheorem{remark}[theorem]{Remark}

\theoremstyle{definition}

\newcommand\bp{\begin{proof}}
\newcommand\ep{\end{proof}}

\newcommand\Dhat{{\hat\Delta}}

\newcommand\A{{\mathcal A}}
\newcommand\CC{{\mathcal C}}

\newcommand\QQ{{\mathcal Q}}
\newcommand\TT{{\mathcal T}}
\newcommand\U{{\mathcal U}}

\newcommand\bb{{\mathfrak b}}
\newcommand\g{{\mathfrak g}}
\newcommand\h{{\mathfrak h}}

\newcommand\sltwo{\mathfrak{sl}_2}

\newcommand\Vect{\mathcal Vec}

\newcommand{\ad}{\operatorname{ad}}

\newcommand{\C}{{\mathbb C}}
\newcommand{\DD}{{\mathbb D}}
\newcommand{\N}{{\mathbb N}}

\newcommand{\R}{{\mathbb R}}

\newcommand\T{{\mathbb T}}
\newcommand\Z{{\mathbb Z}}

\newcommand\eps{\varepsilon}

\newcommand\enu[1]{\smallskip\newline\makebox[5mm][l]{\rm(#1)}}

\begin{document}

\title[Quantized algebras of functions]{Quantized algebras of functions on homogeneous spaces with Poisson stabilizers}

\author[S. Neshveyev]{Sergey Neshveyev}
\address{Department of Mathematics, University of Oslo,
P.O. Box 1053 Blindern, NO-0316 Oslo, Norway}

\email{sergeyn@math.uio.no}

\author[L. Tuset]{Lars Tuset}
\address{Faculty of Engineering, Oslo University College,
P.O. Box 4 St.~Olavs plass, NO-0130 Oslo, Norway}
\email{Lars.Tuset@iu.hio.no}

\thanks{Supported by the Research Council of Norway.}

\date{March 22, 2011; minor corrections April 15, 2011}

\begin{abstract}
Let $G$ be a simply connected semisimple compact Lie group with standard Poisson structure, $K$ a closed Poisson-Lie subgroup, $0<q<1$. We study a quantization $C(G_q/K_q)$ of the algebra of continuous functions on $G/K$. Using results of Soibelman and Dijkhuizen-Stokman we classify the irreducible representations of $C(G_q/K_q)$ and obtain a composition series for $C(G_q/K_q)$. We describe closures of the symplectic leaves of $G/K$ refining the well-known description in the case of flag manifolds in terms of the Bruhat order. We then show that the same rules describe the topology on the spectrum of $C(G_q/K_q)$. Next we show that the family of C$^*$-algebras $C(G_q/K_q)$, $0<q\le1$, has a canonical structure of a continuous field of C$^*$-algebras and provides a strict deformation quantization of the Poisson algebra $\C[G/K]$. Finally, extending a result of Nagy, we show that $C(G_q/K_q)$ is canonically KK-equivalent to $C(G/K)$.
\end{abstract}

\maketitle

\bigskip


\section*{Introduction}

Following the foundational works of Woronowicz \cite{W} and Soibelman and Vaksman \cite{SV1}, the algebras of functions on $q$-deformations of compact groups and their homogeneous spaces were extensively studied in the 90s. Later the interest moved more towards noncommutative geometry of these quantum spaces, see for example \cite{CS}, \cite{DAD}, \cite{NT2} and references therein, leaving the basic algebraic results scattered in the literature and proved mostly in particular cases with various degrees of generality limited often to $SU(2)$ or $SU(N)$ and some of their homogeneous spaces, and more rarely to classical compact simple groups and the corresponding full flag manifolds. The goal of this paper is to establish the main properties of quantized algebras of functions in full generality, for arbitrary Poisson homogeneous spaces of compact semisimple Lie groups such that the stabilizer of one point is a Poisson-Lie subgroup. (Full generality is of course a relative term here, as such spaces form only a relatively small class within the class of all Poisson homogeneous spaces \cite{Ka}.) As often happens, working in the general setting streamlines arguments and renders proofs more transparent, since there is less motivation and possibilities to make use of particular generators and relations. Most results are achieved by first considering $SU(2)$ and then using an inductive argument on the length of an element of the Weyl group. As such the proofs owe much to the fundamental work of Soibelman~\cite{So}. We will now describe the contents of the paper.

In Section~\ref{s1} we briefly remind how a standard Poisson structure on a compact simply connected semisimple Lie group $G$ is defined and what the symplectic leaves are for this structure~\cite{LW,So}. Here we also classify all closed Poisson-Lie subgroups of $G$, a result which probably is known to experts on Poisson geometry.

In Section~\ref{s2} we fix a closed Poisson-Lie subgroup $K$ of $G$ and define the C$^*$-algebra $C(G_q/K_q)$ of functions on the $q$-deformation of $G/K$. The irreducible representations of $C(G_q)$ were classified by Soibelman~\cite{So}. Using his results Dijkhuizen and Stokman~\cite{DS}, following an earlier work of Podkolzin and Vainerman~\cite{PV} on quantum Stiefel manifolds, classified the irreducible representations of quantized algebras of functions on flag manifolds. From this we easily obtain a classification of the irreducible representations of $C(G_q/K_q)$, showing in particular that the equivalence classes of the irreducible representations are in a one-to-one correspondence with the symplectic leaves of $G/K$.

The structure of irreducible representations is refined in Section~\ref{s3}, where we obtain a composition series for $C(G_q/K_q)$. Such a composition series appeared already in work of Soibelman and Vaksman~\cite{SV2} on quantum odd-dimensional spheres. Similar results were then obtained in a number of particular cases~\cite{KL,Sh}, most recently for quantum Stiefel manifolds~\cite{CS}. The main part of the proof can be thought of as an analogue of the fact that the product of symplectic leaves of dimensions~$n$ and~$m$ in~$G$ decomposes into leaves of dimensions $\le n+m$.

Further refinement is obtained in Section~\ref{s4}, where we describe the Jacobson topology on the spectrum of $C(G_q/K_q)$. For $C(G_q)$, when the spectrum is identified with $W\times T$, where $W$ is the Weyl group and $T$ is the maximal torus in $G$, it was observed already by Soibelman~\cite{So} that the closure of $\{w\}\times T$ coincides with the set $\{\sigma\mid \sigma\le w\}\times T$, where $W$ is given the Bruhat order. It follows that the closure of a point $(w,t)\in W\times T$ is a union of sets $\{\sigma\}\times tT_{\sigma,w}$ with $\sigma\le w$ and $T_{\sigma,w}\subset T$. In Section~\ref{s4} we give a combinatorial description of the sets $T_{\sigma,w}$. The corresponding result for $q=1$ is that the closure of the symplectic leaf $\Sigma_w$ associated with $w\in W$ is the union of the sets $\Sigma_\sigma T_{\sigma,w}$ with $\sigma\le w$. This refines the well-known description of the cellular decomposition of $G/T$~\cite{St}.

In the formal deformation setting the algebra $\C[G_q]$, $q=e^{-h}$, of regular functions on $G_q$ is a deformation quantization of the Poisson algebra $\C[G]$. An accepted analytic analogue of deformation quantization is Rieffel's notion of strict deformation quantization~\cite{Ri}. In Section~\ref{s5} we show that the family of C$^*$-algebras $C(G_q/K_q)$ has a canonical continuous field structure, and then that $\C[G_q/K_q]$ define (non-canonically) a strict deformation quantization of $\C[G/K]$. This was proved for $G=SU(2)$ in~\cite{Sh0} and~\cite{Bau} and for $G=SU(N)$ in~\cite{Na3} (although it is not clear from the argument in~\cite{Na3} which Poisson structure on $SU(N)$ is quantized). The main observation which allows us to reduce the proof to the case $G=SU(2)$, and which we already essentially made in~\cite{NT6}, is that it is possible to canonically define a $C[a,b]$-algebra $\Gamma_{alg}((\C[G_q])_{q\in[a,b]})$ playing the role of $\C[G_q]$ when $q$ is considered not as a fixed number, but as the identity function on $[a,b]$. Furthermore, these algebras have the expected functorial properties: given an embedding $K\hookrightarrow G$ we get a homomorphism $\Gamma_{alg}((\C[G_q])_{q})\to \Gamma_{alg}((\C[K_q])_{q})$.

A composition series similar to the one obtained in Section~\ref{s3} was used already by Soibelman and Vaksman~\cite{SV2} to compute the K-theory of the odd-dimensional quantum spheres. Such series were later used for K-theoretic computations in~\cite{Sh} and~\cite{CS}. The most powerful result of this sort was obtained by Nagy~\cite{Na}, who showed that the C$^*$-algebra $C(SU_q(N))$ is KK-equivalent to $C(SU(N))$, and remarked that similar arguments work for all classical simple compact groups. In Section~\ref{s6} we extend this result by showing that $C(G_q/K_q)$ is canonically KK-equivalent to $C(G/K)$.

\bigskip


\section{Poisson-Lie subgroups} \label{s1}

Let $G$ be a simply connected semisimple compact Lie group, $\g$ its
complexified Lie algebra. The universal enveloping algebra $U\g$ is a Hopf $*$-algebra with involution such that the real Lie algebra of~$G$ consists of skew-adjoint elements. Fix a nondegenerate symmetric $\ad$-invariant form on $\g$ such that its restriction to the real Lie algebra of $G$ is negative definite. Let $\h\subset\g$ be the Cartan subalgebra defined by a maximal torus $T$ in $G$. For every root $\alpha\in\Delta$ put $d_\alpha=(\alpha,\alpha)/2$. Let $H_\alpha\in\h$ be the element corresponding to the coroot $\alpha^\vee=2\alpha/(\alpha,\alpha)$ under the identification $\h\cong\h^*$. Under the same identification let $h_\beta\in\h$ be the element corresponding to $\beta\in\h^*$, so $h_\alpha=d_\alpha H_\alpha$ for $\alpha\in\Delta$. Fix a system $\Pi=\{\alpha_1,\dots,\alpha_r\}$ of simple roots. For every positive root $\alpha\in\Delta_+$ choose $E_\alpha\in\g_\alpha$ such that $(E_\alpha,E_\alpha^*)=d_\alpha^{-1}$, and put $F_\alpha=E_\alpha^*\in \g_{-\alpha}$, so that $[E_\alpha,F_\alpha]=H_\alpha$. We write $E_i,F_i,H_i,h_i$ for $E_{\alpha_i},F_{\alpha_i}, H_{\alpha_i}, h_{\alpha_i}$, respectively.  Denote by $\omega_1,\dots,\omega_r$ the fundamental weights, so $\omega_i(H_j)=\delta_{ij}$.

The standard Poisson structure on $G$ is defined by the classical $r$-matrix
$$
r=i\sum_{\alpha\in\Delta_+}d_\alpha(F_\alpha\otimes E_\alpha-E_\alpha\otimes F_\alpha),
$$
meaning that if we consider the Hopf $*$-algebra $\C[G]$ of regular functions on $G$ as a subspace of~$(U\g)^*$, the Poisson bracket on $G$ is given by
\begin{equation} \label{epoissonbracket}
\{f_1,f_2\}=(f_1\otimes f_2)([\Dhat(\cdot),r]) \ \ \hbox{for}\ \ f_1,f_2\in\C[G],
\end{equation}
where $\Dhat$ is the comultiplication on $U\g$.

Soibelman~\cite{So} gave the following description of the symplectic leaves of $G$. For every simple root~$\alpha$ consider the corresponding embedding $\gamma_\alpha\colon SU(2)\to G$. It is a Poisson map when~$SU(2)$ is equipped with the Poisson structure defined by the classical $r$-matrix $id_\alpha(F\otimes E-E\otimes F)$. Consider the symplectic leaf
$$
\Sigma_0=\left\{\begin{pmatrix}\bar z & (1-|z|^2)^{1/2}\\ -(1-|z|^2)^{1/2} & z\end{pmatrix}: |z|<1\right\}\subset SU(2).
$$
Let $W$ be the Weyl group of $G$. Denote by $s_\gamma\in W$ the reflection defined by $\gamma\in\Delta$. We write $s_i$ for~$s_{\alpha_i}$. For every $w\in W$ choose a reduced expression $w=s_{i_1}\dots s_{i_n}$ and consider the map
\begin{equation} \label{esymplleaf}
\gamma_w\colon\Sigma_0^{n}\to G, \ \ (g_1,\dots,g_n)\mapsto \gamma_{i_1}(g_1)\dots \gamma_{i_n}(g_n),
\end{equation}
where $\gamma_i=\gamma_{\alpha_i}$; for $w=e$ the image of $\gamma_e$ consists solely of the identity element in $G$. It is a symplectomorphism of~$\Sigma_0^{n}$ onto a symplectic leaf $\Sigma_w$ of $G$. The leaf $\Sigma_w$ does not depend on the reduced expression for~$w$, although the map $\gamma_w$ depends on it. The decomposition of $G$ into its symplectic leaves is given by $G=\sqcup_{w\in W,t\in T}\Sigma_wt$.

We next define a class of subgroups of $G$. Let $S$ be a subset of $\Pi$. Denote by $\tilde K^S$ the closed connected subgroup of $G$ such that its complexified Lie algebra $\tilde\g_S$ is generated by the elements $E_i$ and $F_i$ with $\alpha_i\in S$, so
$$
\tilde\g_S=\operatorname{span}\{H_i\mid \alpha_i\in S\}\oplus\bigoplus_{\alpha\in\Delta_S}\g_\alpha,
$$
where $\Delta_S$ is the set of roots that lie in the group generated by $\alpha_i\in S$. Denote by $P(S^c)$ the subgroup of the weight lattice $P$ generated by the fundamental weights $\omega_i$ with $\alpha_i\in S^c=\Pi\setminus S$. Let $L$ be a subgroup of $P(S^c)$. Identifying $P$ with the dual group of the maximal torus $T$, denote by $T_L$ the annihilator of $L$ in $T$. Since $T$ normalizes $\tilde K^S$, the group $K^{S,L}$ generated by $\tilde K^S$ and $T_L$ is a closed subgroup of $G$, and its complexified Lie algebra is
$$
\g_{S,L}=\h_L\oplus\bigoplus_{\alpha\in\Delta_S}\g_\alpha,
$$
where $\h_L\subset\h$ is the annihilator of $L\subset\h^*$. Note that if $L=P(S^c)$ then $K^{S,L}$ is the group $\tilde K^S$. If $L=0$, we write $K^S$ for $K^{S,L}$. Then $K^S=G\cap P_S$, where $P_S\subset G_\C$ is the parabolic subgroup corresponding to $S$, and $\tilde K^S$ is the semisimple part of $K^S$.

\begin{proposition} \label{ppoissonsub}
For any subset $S\subset\Pi$ and any subgroup $L\subset P(S^c)$ we have:
\enu{i} $K^{S,L}$ is a Poisson-Lie subgroup of $G$, and any closed Poisson-Lie subgroup of $G$ is of this form for uniquely defined $S$ and $L$;
\enu{ii} $K^{S,L}\cap T=T_L$ and $(K^{S,L})^\circ\cap T=T_L^\circ$;
\enu{iii} $K^{S,L}$ is connected if and only if $P(S^c)/L$ is torsion-free.
\end{proposition}

\bp (i) By \cite[Proposition 2.1]{Stok} a closed connected Lie subgroup $K$ of $G$ is a Poisson-Lie subgroup if and only if its complexified Lie algebra $\mathfrak k_\C$ lies between $\tilde \g_S$ and $\g_S$ for some $S\subset \Pi$, so it has the form
\begin{equation} \label{epoissonsub}
\mathfrak k_\C=\mathfrak a\oplus\bigoplus_{\alpha\in\Delta_S}\g_\alpha,
\end{equation}
where $\mathfrak a$ is the complexified Lie algebra of $K\cap T$. It follows that for any $S\subset\Pi$ and $L\subset P(S^c)$, $(K^{S,L})^\circ$ is a Poisson-Lie subgroup. Furthermore, by construction $K^{S,L}$ is a finite disjoint union of sets of the form $(K^{S,L})^\circ t$ with $t\in T_L$. Since the translations by the elements of $T$ are Poisson maps, every such set $(K^{S,L})^\circ t$ is a Poisson submanifold of $G$, hence $K^{S,L}$ is a Poisson-Lie subgroup.

Conversely, assume $K$ is a closed Poisson-Lie subgroup of $G$. Assume first that $K$ is connected. Then its complexified Lie algebra has form \eqref{epoissonsub}. Denote by $L$ the annihilator of $K\cap T$ in $\hat T=P$. Then $K\cap T=T_L$ and $\mathfrak a=\h_L$, and since $H_i$ lies in this Lie algebra for $\alpha_i\in S$, we have $L\subset P(S^c)$. Therefore $K=(K^{S,L})^\circ$. Observe next that the group $T_L=K\cap T$ is connected, since it is abelian and contains a maximal torus of $K$. Hence the group $K^{S,L}=\tilde K^S T_L$ is connected, and thus $K=K^{S,L}$.

Consider now a not necessarily connected closed Poisson-Lie subgroup $K$. Then $K^\circ=K^{S,\Gamma}$ for some  $S\subset\Pi$ and $\Gamma\subset P(S^c)$. Let $g\in K$. Consider a symplectic leaf $\Sigma$ of $G$ passing through $g$. By assumption the whole leaf $\Sigma$, and hence its closure, lies in $K$. From the description of symplectic leaves given above it is clear that $\bar\Sigma\cap T\ne\emptyset$. Since $\bar\Sigma$ is connected, it follows that there exists $t\in K\cap T$ such that $gt^{-1}\in K^\circ$. Therefore $K$ is generated by $K^\circ=K^{S,\Gamma}$ and $K\cap T$. Let $L\subset \Gamma$ be such that $K\cap T=T_L$. Then we conclude that $K=K^{S,L}$.

That $S$ is uniquely defined by $K$, is clear. That $L$ is also uniquely defined, will follow from (ii).

\smallskip

(ii) As we already observed, the group $(K^{S,L})^\circ\cap T$ is connected. Since the Lie algebras of $K^{S,L}\cap T$ and $T_L$ clearly coincide, we conclude that $(K^{S,L})^\circ\cap T=T_L^\circ$. Since $K^{S,L}=(K^{S,L})^\circ T_L$, it follows that $K^{S,L}\cap T=T_L$.

\smallskip

(iii) As $K^{S,L}=(K^{S,L})^\circ T_L$ and $(K^{S,L})^\circ\cap T=T_L^\circ$, we have $K^{S,L}/(K^{S,L})^\circ=T_L/T_L^\circ$. In particular, $K^{S,L}$ is connected if and only if $T_L$ is connected, that is, $\hat T_L\cong P/L$ is torsion-free, or equivalently, $P(S^c)/L$ is torsion-free.
\ep

To describe the symplectic leaves of the Poisson manifold $G/K^{S,L}$, consider the subgroup  $W_S$ of~$W$ generated by the simple reflections $s_\alpha$ with $\alpha\in S$. Let $W^S\subset W$ be the set of elements $w$ such that $w(\alpha)>0$ for all $\alpha\in S$. Then (see e.g.~page~140 in \cite{St}) every element $w\in W$ decomposes uniquely as $w=w'w''$ with $w'\in W^S$ and $w''\in W_S$, and we have $\ell(w)=\ell(w')+\ell(w'')$; recall also that the length of an element in~$W_S$ is the same as in $W$.

\begin{proposition}
Let $\pi\colon G\to G/K^{S,L}$ be the quotient map. Then, for every $w\in W^S$ and $t\in T$, the map~$\pi$ defines a symplectomorphism of $\Sigma_w t$ onto a symplectic leaf of $G/K^{S,L}$. The decomposition of~$G/K^{S,L}$ into its symplectic leaves is given by $\sqcup_{w\in W^S,t\in T/T_L}\pi(\Sigma_wt)$.
\end{proposition}

\bp This is just a slight extension of results of Lu and Weinstein~\cite{LW} and Soibelman~\cite{So}, see also~\cite{DS}. We have the decompositions $G=\sqcup_{w\in W,t\in T}\Sigma_wt$ and $K^{S,L}=\sqcup_{w\in W_S,t\in T_L}\Sigma_wt$. Using that $T\Sigma_w=\Sigma_wT$ for any $w\in W$, and that the multiplication map $\Sigma_{w'}\times\Sigma_{w''}\to\Sigma_{w'w''}$ is a bijection for $w'\in W^S$ and $w''\in W_S$ (since $\ell(w'w'')=\ell(w')+\ell(w'')$), we conclude that $\pi$ is injective on every leaf $\Sigma_w t$ with $w\in W^S$ and arbitrary $t\in T$, and $G/K^{S,L}=\sqcup_{w\in W^S,t\in T/T_L}\pi(\Sigma_wt)$. Since by \cite[Theorem~4.6]{LW} the symplectic leaves of $G$ and $G/K^{S,L}$ are orbits of the right dressing action of the Poisson-Lie dual of $G$, the sets $\pi(\Sigma_wt)$ are symplectic leaves of $G/K^{S,L}$ for all $w\in W$ and $t\in T$.
\ep

\bigskip


\section{Irreducible representations of quantized function algebras} \label{s2}

Fix $q\in(0,1]$. If $q=1$ we put $U_1\g=U\g$. For $q\ne1$ the quantized universal
enveloping algebra~$U_q\g$ is generated by elements $E_i$, $F_i$, $K_i$,
$K_i^{-1}$, $1\le i\le r$, satisfying the relations
$$
K_iK_i^{-1}=K_i^{-1}K_i=1,\ \ K_iK_j=K_jK_i,\ \
K_iE_jK_i^{-1}=q_i^{a_{ij}}E_j,\ \
K_iF_jK_i^{-1}=q_i^{-a_{ij}}F_j,
$$
$$
E_iF_j-F_jE_i=\delta_{ij}\frac{K_i-K_i^{-1}}{q_i-q_i^{-1}},
$$
$$
\sum^{1-a_{ij}}_{k=0}(-1)^k\begin{bmatrix}1-a_{ij}\\
k\end{bmatrix}_{q_i} E^k_iE_jE^{1-a_{ij}-k}_i=0,\ \
\sum^{1-a_{ij}}_{k=0}(-1)^k\begin{bmatrix}1-a_{ij}\\
k\end{bmatrix}_{q_i} F^k_iF_jF^{1-a_{ij}-k}_i=0,
$$
where $\displaystyle\begin{bmatrix}m\\
k\end{bmatrix}_{q_i}=\frac{[m]_{q_i}!}{[k]_{q_i}![m-k]_{q_i}!}$,
$[m]_{q_i}!=[m]_{q_i}[m-1]_{q_i}\dots [1]_{q_i}$,
$\displaystyle[n]_{q_i}=\frac{q_i^n-q_i^{-n}}{q_i-q_i^{-1}}$,
$q_i=q^{d_i}$ and $d_i=d_{\alpha_i}$. This is a Hopf $*$-algebra with coproduct $\Dhat_q$ and
counit $\hat\eps_q$ defined by
$$
\Dhat_q(K_i)=K_i\otimes K_i,\ \
\Dhat_q(E_i)=E_i\otimes1+ K_i\otimes E_i,\ \
\Dhat_q(F_i)=F_i\otimes K_i^{-1}+1\otimes F_i,
$$
$$
\hat\eps_q(E_i)=\hat\eps_q(F_i)=0,\ \ \hat\eps_q(K_i)=1,
$$
and with involution given by $K_i^*=K_i$, $E_i^*=F_iK_i$, $F_i^*=K_i^{-1}E_i$.

If $V$ is a finite dimensional $U_q\g$-module and $\lambda\in
P$ is an integral weight, denote by $V(\lambda)$ the
space of vectors $v\in V$ of weight $\lambda$, so that
$K_iv=q^{(\lambda,\alpha_i)}v=q_i^{(\lambda,\alpha_i^\vee)}v$ for all $i$. Recall that $V$ is called admissible if
$V=\oplus_{\lambda\in P}V(\lambda)$. We denote by $\CC_q(\g)$ the tensor category of finite dimensional admissible $U_q\g$-modules.

Denote by $\C[G_q]\subset (U_q\g)^*$ the Hopf $*$-algebra of matrix coefficients of finite dimensional admissible $U_q\g$-modules, and let $C(G_q)$ be its C$^*$-enveloping algebra.

Consider also the endomorphism ring $\U(G_q)$ of the forgetful functor $\CC_q(\g)\to\Vect$.
In other words, if for every $\lambda\in P_+$ we fix
an irreducible $*$-representation of $U_q\g$ on a Hilbert space $V_\lambda$ with highest weight
$\lambda$, then  $\U(G_q)$ can be identified with $\prod_{\lambda\in P_+}B(V_\lambda)$. Yet another way to think of~$\U(G_q)$ is as the algebra of closed densely defined operators affiliated with the von Neumann algebra~$W^*(G_q)$ of $G_q$. The maximal torus $T\subset G$ can be considered as a subset of group-like elements of $W^*(G_q)\subset \U(G_q)$: if $X\in\mathfrak t$ then for any admissible $U_q\g$-module $V$ and $\lambda\in P$ the element $\exp( X)\in T$ acts on $V(\lambda)$ as multiplication by $e^{\lambda(X)}$. Under this embedding $T\hookrightarrow\U(G_q)$, we have $K_j^{it}=\exp(it(\log q) h_j)\in T$ for $j=1,\dots,r$ and $t\in\R$.

From now on fix a subset $S\subset\Pi$ and a subgroup $L\subset P(S^c)$. Let $\U(K^{S,L}_q)$ be the $\sigma(\U(G_q),\C[G_q])$-closed subalgebra of $\U(G_q)$ generated by $T_L$ and $E_i$, $F_i$ with $\alpha_i\in S$. In other words, an element $\omega\in\U(G_q)$ belongs to $\U(K^{S,L}_q)$ if and only if for every finite dimensional admissible $U_q\g$-module $V$ the operator of the action by $\omega$ on $V$ lies in the algebra generated by $T_L$ and $E_i$, $F_i$ with $\alpha_i\in S$. Denote by $\C[K^{S,L}_q]\subset \U(K^{S,L}_q)^*$ the Hopf $*$-algebra that is the image of $\C[G_q]$ under the restriction map $\U(G_q)^*\to\U(K^{S,L}_q)^*$, and let $C(K^{S,L}_q)$ be its C$^*$-enveloping algebra. By construction we have an epimorphism $\pi\colon\C[G_q]\to\C[K^{S,L}_q]$ of Hopf $*$-algebras. Put
$$
\C[G_q/K^{S,L}_q]=\{a\in \C[G_q]\mid (\iota\otimes\pi)\Delta_q(a)=a\otimes 1\},
$$
where $\Delta_q$ is the comultiplication on $\C[G_q]$. Equivalently, $a\in \C[G_q/K^{S,L}_q]$ if and only if $(\iota\otimes \omega)\Delta_q(a)=\hat\eps_q(\omega)a$ for all $\omega\in\U(K^{S,L}_q)$. Denote by $C(G_q/K^{S,L}_q)$ the norm-closure of  $\C[G_q/K^{S,L}_q]$ in $C(G_q)$.

For $\lambda\in P_+$ and $\xi,\zeta\in V_\lambda$ denote by $C^\lambda_{\zeta,\xi}\in \C[G_q]$ the matrix coefficient $(\cdot \,\xi,\zeta)$. Then $\C[G_q/K^{S,L}_q]$ is the linear span of elements $C^\lambda_{\zeta,\xi}$ such that $\lambda\in P_+$, $\zeta\in V_\lambda$ and $\xi\in V_\lambda$ is fixed by $K^{S,L}_q$, that is, $\omega\xi=\hat\eps_q(\omega)\xi$ for all $\omega\in\U(K^{S,L}_q)$.

\begin{remark}
The above description of $\C[G_q/K^{S,L}_q]$ as the linear span of certain elements $C^\lambda_{\zeta,\xi}$ implies that there exists a C$^*$-enveloping algebra of $\C[G_q/K^{S,L}_q]$, since if $\{ e_i\}_i$ is an orthonormal basis in $V_\lambda$ then $\sum_i(C^\lambda_{e_i,\xi})^*C^\lambda_{e_i,\xi}=\|\xi\|^21$, so that the norm of $C^\lambda_{\zeta,\xi}\in\C[G_q/K^{S,L}_q]$ in every representation of $\C[G_q/K^{S,L}_q]$ by bounded operators is not bigger than $\|\zeta\|\,\|\xi\|$. This C$^*$-algebra coincides with~$C(G_q/K^{S,L}_q)$. This was proved by Stokman~\cite{Stok} in the case $L=0$ using results of~\cite{DS} on representations of a certain subalgebra of $\C[G_q/K^{S,L}_q]$. An alternative way to see this is to use coamenability of $G_q$~\cite{Ba1}, see also Appendix~A in~\cite{NT5}, together with the following well-known fact: if a coamenable compact quantum group $G$ acts ergodically on a unital C$^*$-algebra $A$ (that is, we have a coaction $\alpha\colon A\to C(G)\otimes A$ such that $A^G=\{a\in A\mid\alpha(a)=1\otimes a\}=\C1$), $\A\subset A$ is the $*$-subalgebra spanned by the spectral subspaces of the action, and there exists an enveloping C$^*$-algebra $\tilde A$ for~$\A$, then~$\tilde A=A$. Indeed, let $\tilde\alpha\colon \tilde A\to C(G)\otimes \tilde A$ be the action of $G$ extending that on~$\A$ and let $\pi\colon\tilde A\to A$ be the quotient map. Consider the conditional expectation $\tilde E\colon \tilde A\to \tilde A^G=\C1$ defined by $\tilde E(a)=(h\otimes\iota)\tilde\alpha(a)$, where $h$ is the Haar state on $C(G)$, and a similar conditional expectation $E\colon A\to A^G=\C1$. As $h$ is faithful by coamenability of $G$, these conditional expectations are faithful. Since $\pi\tilde E=E\pi$, it follows that the kernel of $\pi$ is trivial.
\end{remark}

Our goal is to describe the irreducible representations of the C$^*$-algebra $C(G_q/K^{S,L}_q)$ for $q\in(0,1)$. For this recall the classification of irreducible representations of $C(G_q)$ obtained by Soibelman~\cite{So}.

Consider first the case $G=SU(2)$. We assume that the invariant symmetric form on $\sltwo(\C)$ is the standard one, so $(\alpha,\alpha)=2$ for the unique simple root $\alpha$. Consider the fundamental representation
$$
E\mapsto\begin{pmatrix}0 & q^{1/2}\\ 0 & 0\end{pmatrix}, \ \
F\mapsto\begin{pmatrix}0 & 0\\ q^{-1/2} & 0\end{pmatrix}, \ \
K\mapsto\begin{pmatrix}q & 0\\ 0 & q^{-1}\end{pmatrix}
$$
of $U_q\sltwo$. Then the corresponding corepresentation of $\C[SU_q(2)]$ has the form
$$
\begin{pmatrix}\alpha & -q\gamma^*\\ \gamma & \alpha^*\end{pmatrix},
$$
and the elements $\alpha,\gamma\in\C[SU_q(2)]$ satisfy the relations
$$
\alpha^*\alpha+\gamma^*\gamma=1,\ \, \alpha\alpha^*+q^2\gamma^*\gamma=1, \ \ \gamma^*\gamma=\gamma\gamma^*,
\ \ \alpha\gamma=q\gamma\alpha, \ \ \alpha\gamma^*=q\gamma^*\alpha.
$$
We will write $\alpha_q,\gamma_q$ when we want to emphasize that these are elements of $\C[SU_q(2)]$ for a particular~$q$.

Define a representation $\rho_q$ of $C(SU_q(2))$ on $\ell^2(\Z_+)$ by
\begin{equation} \label{erepSU2}
\rho_q(\alpha)e_n=\sqrt{1-q^{2n}}\,e_{n-1}, \ \ \rho_q(\gamma)e_n=-q^ne_n, \ \ n\ge0.
\end{equation}

Return to the general case. For every $1\le i\le r$ consider the homomorphism $\sigma_i\colon C(G_q)\to C(SU_{q_i}(2))$ which is dual to the embedding $U_{q_i}\sltwo\hookrightarrow U_q\g$ corresponding to the simple root $\alpha_i$. Then $\pi_i=\rho_{q_i}\sigma_i$ is a representation of $C(G_q)$ on $\ell^2(\Z_+)$. Now for every element $w\in W$ fix a reduced decomposition $w=s_{i_1}\dots s_{i_n}$ and put
$$
\pi_w=\pi_{i_1}\otimes\dots\otimes\pi_{i_n},
$$
so $\pi_w(a)=(\pi_{i_1}\otimes\dots\otimes\pi_{i_n})\Delta_q^{(n-1)}(a)$. Then $\pi_w$ is a representation of $C(G_q)$ on $\ell^2(\Z_+)^{\otimes\ell(w)}=\ell^2(\Z_+^{\ell(w)})$. Up to equivalence it does not depend on the choice of the reduced expression for $w$. In addition we have one-dimensional representations $\pi_t$ of $C(G_q)$ defined by the points of the maximal torus $T\subset\U(G_q)=\C[G_q]^*$. In other words, $\pi_t(C^\lambda_{\zeta,\xi})=(t\xi,\zeta)$. Then the result of Soibelman says that the representations $\pi_{w,t}=\pi_w\otimes\pi_t$ are irreducible, mutually inequivalent, and exhaust all irreducible representations of $C(G_q)$ up to equivalence. Note that
\begin{equation}\label{et}
\pi_{w,t}(C^\lambda_{\zeta,\xi})=\pi_{w}(C^\lambda_{\zeta,t\xi}).
\end{equation}

The following result is a minor generalization of \cite[Theorem~5.9]{DS}.

\begin{theorem} \label{treps}
Assume $q\in(0,1)$. Then for every $w\in W^S$ and $t\in T$ the restriction of the representation $\pi_{w,t}$ of $C(G_q)$ to $C(G_q/K^{S,L}_q)$ is irreducible. Such representations exhaust all irreducible representations of $C(G_q/K^{S,L}_q)$ up to equivalence.
For $w,w'\in W^S$ and $t,t'\in T$ the restrictions of~$\pi_{w,t}$ and~$\pi_{w',t'}$ to $C(G_q/K^{S,L}_q)$ are equivalent if and only if $w=w'$ and $t't^{-1}\in T_L$, and in this case they are actually equal.
\end{theorem}

To prove the theorem we will need further properties of the representations $\pi_w$. Let $\lambda\in P_+$.
Fix a highest weight unit vector $\xi_\lambda\in V_\lambda$. For every $w\in W$ choose a unit vector $\eta\in V_\lambda$ of weight $w\lambda$. Since the weight spaces $V_\lambda(\lambda)$ and $V_\lambda(w\lambda)$ are one-dimensional, the element $C^\lambda_{\eta,\xi_\lambda}$ does not depend on the choice of $\xi_\lambda$ and $\eta$ up to a factor of modulus one. To simplify the notation we will thus write~$C^\lambda_{w\lambda,\lambda}$ for~$C^\lambda_{\eta,\xi_\lambda}$.

\begin{lemma}[\cite{So}] \label{lsoib}
Let $w\in W$ and $\lambda\in P_+$. Then
\enu{i} $\pi_w(C^\lambda_{w\lambda,\lambda})$ is a compact contractive diagonalizable operator with zero kernel, and the vector $e_0^{\otimes\ell(w)}\in\ell^2(\Z_+)^{\otimes\ell(w)}$ is its only (up to a scalar factor) eigenvector with eigenvalue of modulus $1$;
\enu{ii} if $\zeta\in V_\lambda$ is orthogonal to $(U_q\bb) V_\lambda(w\lambda)$, where $U_q\bb\subset U_q\g$ is the subalgebra generated by $K_i$, $K_i^{-1}$ and $E_i$, $1\le i\le r$, then $\pi_w(C^\lambda_{\zeta,\xi_\lambda})=0$.
\end{lemma}

\bp Part (i) is a consequence of the proof of \cite[Proposition~6.1.5]{KS}, see also identity (6.2.4) there (although notice that a factor of modulus one depending on the choice of orthonormal bases is missing there).


Part (ii) is \cite[Theorem~6.2.1]{KS}, since by that theorem $\pi_w$ corresponds to the Schubert cell $X_w$ in the terminology of~\cite{KS}, which in particular means that it has the property in the statement of the lemma.
\ep

\bp[Proof of Theorem~\ref{treps}]
Write $K^S_q$ for $K^{S,0}_q$. By \cite[Theorem~5.9]{DS} the restrictions of the representations~$\pi_{w,t}$ to $C(G_q/K^S_q)\subset C(G_q/K_q^{S,L})$ are irreducible for $w\in W^S$. Hence the restrictions of $\pi_{w,t}$ to~$C(G_q/K_q^{S,L})$ are irreducible as well. To see that this way we get all irreducible representations, note that any irreducible representation of $C(G_q/K_q^{S,L})$ extends to an irreducible representation of~$C(G_q)$ on a larger space. Therefore we have to find decompositions of $\pi_{w,t}$ into irreducible representations of $C(G_q/K_q^{S,L})$ for arbitrary $w\in W$ and $t\in T$. Write $w=w'w''$ with $w'\in W^S$ and $w''\in W_S$. We may assume that $\pi_w=\pi_{w'}\otimes\pi_{w''}$. Then by the proof~\cite[Proposition~5.7]{DS} we have
\begin{equation} \label{eDS}
\pi_w(a)=\pi_{w'}(a)\otimes 1^{\otimes\ell(w'')}\ \ \hbox{for all}\ \ a\in C(G_q/K_q^{S,L})
\end{equation}
(the key point here is that for the homomorphism $\sigma_i\colon C(G_q)\to C(SU_{q_i}(2))$ we have $\sigma_i(a)=\eps_q(a)1$ for all $a\in C(G_q/K_q^{S,L})$ and $\alpha_i\in S$, where $\eps_q$ is the counit on $C(G_q)$). Using \eqref{et} it follows that if $\zeta\in V_\lambda$ and $\xi\in V_\lambda$ is fixed by $K^{S,L}_q$ then
$$
\pi_{w,t}(C^\lambda_{\zeta,\xi})=\pi_{w}(C^\lambda_{\zeta,t\xi})=\pi_{w'}(C^\lambda_{\zeta,t\xi})\otimes 1^{\otimes\ell(w'')}
=\pi_{w',t}(C^\lambda_{\zeta,\xi})\otimes 1^{\otimes\ell(w'')}.
$$
We therefore see that the representations $\pi_{w',t}$ with $w'\in W^S$ and $t\in T$ exhaust all irreducible representations of $C(G_q/K_q^{S,L})$.

Consider now two representations $\pi_{w,t}$ and $\pi_{w',t'}$ with $w,w'\in W^S$ and $t,t'\in T$. By~\cite[Theorem~5.9]{DS} if $w\ne w'$ then already the restrictions of these representations to $C(G_q/K^S_q)\subset C(G_q/K_q^{S,L})$ are inequivalent. Therefore assume $w=w'$. Since $\pi_{w,t}$ and $\pi_{w,t'}$ coincide on $C(G_q/K^S_q)$ and are irreducible as representations of $C(G_q/K^S_q)$, they can be equivalent as representations of $C(G_q/K_q^{S,L})$ only if they coincide on $C(G_q/K_q^{S,L})$. If $t't^{-1}\in T_L$, this is indeed the case, since for $\zeta\in V_\lambda$ and $\xi\in V_\lambda$ fixed by $K^{S,L}_q$ we have
$$
\pi_{w,t}(C^\lambda_{\zeta,\xi})=\pi_{w}(C^\lambda_{\zeta,t\xi})=\pi_{w}(C^\lambda_{\zeta,t'\xi})=
\pi_{w,t'}(C^\lambda_{\zeta,\xi}).
$$
Assume now that $t't^{-1}\notin T_L$. Then there exists $\nu\in L\subset P=\hat T$ such that $\nu(t't^{-1})\ne1$. Choose weights $\lambda,\mu\in P_+(S^c)=P_+\cap P(S^c)$ such that $\lambda-\mu=\nu$. We have $E_i\xi_\lambda=F_i\xi_\lambda=0$ for $\alpha_i\in S$, so that $(\iota\otimes \omega)\Delta_q(C^\lambda_{w\lambda,\lambda})=\hat\eps_q(\omega)C^\lambda_{w\lambda,\lambda}$ for $\omega$ lying in the algebra generated by $E_i$ and $F_i$ with $\alpha_i\in S$. We also have $(\iota\otimes \tau)\Delta_q(C^\lambda_{w\lambda,\lambda})=\lambda(\tau)C^\lambda_{w\lambda,\lambda}$ for $\tau\in T$. It follows that $(C^\mu_{w\mu,\mu})^*C^\lambda_{w\lambda,\lambda}\in\C[G_q/K^{S,L}_q]$. Using \eqref{et} we get
$$
\pi_{w,t}((C^\mu_{w\mu,\mu})^*C^\lambda_{w\lambda,\lambda})
=\overline{\mu(t)}\lambda(t)\pi_{w}((C^\mu_{w\mu,\mu})^*C^\lambda_{w\lambda,\lambda})
=\nu(t)\pi_w((C^\mu_{w\mu,\mu})^*C^\lambda_{w\lambda,\lambda}),
$$
and similarly $\pi_{w,t'}((C^\mu_{w\mu,\mu})^*C^\lambda_{w\lambda,\lambda})
=\nu(t')\pi_w((C^\mu_{w\mu,\mu})^*C^\lambda_{w\lambda,\lambda})$. Since $\pi_w((C^\mu_{w\mu,\mu})^*C^\lambda_{w\lambda,\lambda})\ne0$ by Lemma~\ref{lsoib}(i), we see that $\pi_{w,t}\ne\pi_{w,t'}$ on $C(G_q/K^{S,L}_q)$.
\ep

\begin{corollary}
The $*$-algebra $\C[G_q/K^{S,L}_q]$ is spanned by the elements of the form $(C^\mu_{\zeta,\xi_\mu})^*C^\lambda_{\eta,\xi_\lambda}$, where $\mu,\lambda\in P_+(S^c)$ are such that $\lambda-\mu\in L$, $\zeta\in V_\mu$ and $\eta\in V_\lambda$.
\end{corollary}

\bp This is similar to \cite[Theorem~2.5]{Stok}. That the linear span $\A$ in the formulation forms a $*$-algebra, is proved exactly as~\cite[Lemma~4.3]{DS}. That $\A\subset C[G_q/K^{S,L}_q]$, is checked in the same way as that $(C^\mu_{w\mu,\mu})^*C^\lambda_{w\lambda,\lambda}\in\C[G_q/K^{S,L}_q]$ in the proof of the above theorem. Since $\A$ is invariant with respect to the left coaction $\Delta_q\colon\C[G_q]\to C(G_q)\otimes \C[G_q]$, we have $\bar\A\cap\C[G_q]=\A$. On the other hand, by the proof of the above theorem and by~\cite[Lemma~5.8]{DS} the algebra $\A$ has the property that the restriction of irreducible representations of $C(G_q/K^{S,L}_q)$ to $\A$ are irreducible and inequivalent representations restrict to inequivalent representations. By the Stone-Weierstrass theorem for type~I C$^*$-algebras it follows that $\bar\A=C(G_q/K^{S,L}_q)$. Hence $\A=\C[G_q/K^{S,L}_q]$.
\ep

\bigskip


\section{Composition series} \label{s3}

In this section we will use the classification of irreducible representations of $C(G_q/K^{S,L}_q)$ to construct a composition series for $C(G_q/K^{S,L}_q)$. Since in the subsequent sections it will be important to have such a sequence for all $q$ including $q=1$, it is convenient to look at the irreducible representations in a slightly different way.

For $q\in[0,1)$ denote by $C(\bar\DD_q)$ the universal unital C$^*$-algebra with one generator $Z_q$ such that
$$
1-Z_q^*Z_q=q^2(1-Z_qZ_q^*).
$$
Since the least upper bounds of the spectra of $aa^*$ and $a^*a$ coincide, it is easy to see that in every nonzero representation of the above relation the norm of $Z_q$ is equal to $1$, so $C(\bar\DD_q)$ is well-defined. It follows then, see e.g. Section V in \cite{KL}, that $C(\bar\DD_q)$ is isomorphic to the Toeplitz algebra $\TT\subset B(\ell^2(\Z_+))$ via an isomorphism which maps $Z_q$ into the operator
$$
e_n\mapsto \sqrt{1-q^{2(n+1)}}\, e_{n+1}, \ \ n\ge0.
$$
The inverse homomorphism maps the operator $S\in\TT$ of the shift to the right into $Z_q(Z_q^*Z_q)^{-1/2}$. Under this isomorphism the representation $\rho_q\colon C(SU_q(2))\to B(\ell^2(\Z_+))$ defined by~\eqref{erepSU2} becomes the $*$-homomorphism $C(SU_q(2))\to C(\bar\DD_q)$ given by
\begin{equation} \label{erho}
\rho_q(\alpha_q)=Z_q^*, \ \ \rho_q(\gamma_q)=-(1-Z_qZ_q^*)^{1/2}.
\end{equation}
In this form $\rho_q$ makes sense for $q=1$. Namely, consider the C$^*$-algebra $C(\bar\DD_1)=C(\bar\DD)$ of continuous functions on the closed unit disk, and denote by $Z_1$ its standard generator, $Z_1(z)=z$. Then $\rho_1\colon C(SU(2))\to C(\bar\DD_1)$ defined by the above formula is the homomorphism of restriction of a function on $SU(2)$ to the closure of the symplectic leaf
$$
\Sigma_0=\left\{\begin{pmatrix}\bar z & (1-|z|^2)^{1/2}\\ -(1-|z|^2)^{1/2} & z\end{pmatrix}: |z|<1\right\}\cong\DD.
$$

The representation $\pi_{w,t}$ defined by $w=s_{i_1}\dots s_{i_n}$ now becomes a $*$-homomorphism
$$
\pi_{w,t}\colon C(G_q)\to C(\bar\DD_{q_{i_1}})\otimes\dots\otimes
C(\bar\DD_{q_{i_n}}).
$$
For $q=1$ this is exactly the homomorphism $\gamma^*_w$, where $\gamma_w\colon\DD^n\cong\Sigma_0^n\to G$ is defined by~\eqref{esymplleaf} and extended by continuity to $\bar\DD^n$.

For every $q\in[0,1]$ we have a $*$-homomorphism $C(\bar\DD_q)\to C(\T)$ mapping $Z_q$ into the standard generator of $C(\T)$. Denote by $C_0(\DD_q)$ its kernel. For $q=1$ this is the usual algebra of continuous functions on $\bar\DD$ vanishing on the boundary. For $q\in[0,1)$ this is the ideal of compact operators in~$\TT=C(\bar\DD_q)$.

We can now formulate the main result of this section.

\begin{theorem} \label{tcompseries}
Assume $q\in(0,1]$. Let $w_0$ be the longest element in the Weyl group. Write  $w_0=w_0'w_0''$ with $w'_0\in W^S$ and $w_0''\in W_S$, and put $m_0=\ell(w_0')$. For every $0\le m\le m_0$ denote by $J_m$ the ideal in $C(G_q/K^{S,L}_q)$ consisting of elements $a$ such that $\pi_{w,t}(a)=0$ for all $t\in T$ and $w\in W^S$ with $\ell(w)=m$. Then
$$
0=J_{m_0}\subset J_{m_0-1}\subset\dots\subset J_0\subset J_{-1}=C(G_q/K^{S,L}_q),
$$
and for every $0\le m\le m_0$ we have
$$
J_{m-1}/J_m\cong\bigoplus_{w\in W^S: \ell(w)=m}C(T/T_L;C_0(\DD_{q_{i_1(w)}})\otimes\dots\otimes
C_0(\DD_{q_{i_m(w)}})), \ \ a\mapsto (a_w)_w,
$$
where $a_w(t)=\pi_{w,t}(a)$ and $w=s_{i_1(w)}\dots s_{i_m(w)}$ is the fixed reduced decomposition of $w$ used to define~$\pi_w$.
\end{theorem}

Note that there is no ambiguity in the definition of $a_w$, since by Theorem~\ref{treps} we have $\pi_{w,t}(a)=\pi_{w,t'}(a)$ if $t't^{-1}\in T_L$.

We need some preparation to prove this theorem. Recall some properties of the Bruhat order, see e.g.~\cite{BGG}.
The Bruhat order on $W$ is defined by declaring that $w\le w'$ iff $w$ can be obtained from $w'$ by dropping letters in some (equivalently any) reduced word for~$w'$. Furthermore, then the letters can be dropped in such a way that we get a reduced word for~$w$.

Consider the coset space $X_S=W/W_S$. For $x\in X_S$ define
$$
\ell_S(x)=\min\{\ell(w)\mid x=wW_S\}.
$$
As we know, every coset $x\in X_S$ has a unique representative $w_x\in W^S$, and $w_x$ is the smallest element in $x$ in the Bruhat order; in particular, $\ell_S(x)=\ell(w_x)$. Define an order on $X_S$ by declaring $x\le y$ iff $w_x\le w_y$, and call it again the Bruhat order.

\begin{lemma} \label{lBruhat}
We have:
\enu{i} the factor-map $W\to X_S$ is order-preserving;
\enu{ii} for any $x\in X_S$ and $\alpha\in\Pi$, either $s_\alpha x=x$ or $w_{s_\alpha x}=s_\alpha w_x$.
\end{lemma}

\bp To prove (i), take $u,v\in W$ such that $u\le v$, and put $x=uW_S$, $y=vW_S$. Write $v=w_yw$ with $w\in W_S$. Then
$w_x\le u\le v=w_yw$. Take a reduced word for $w_yw$ which is the concatenation of a reduced word for $w_y$ and a reduced word for $w$ in letters $s_\alpha$ with $\alpha\in S$. Then a reduced word for $w_x$ can be obtained by dropping some letters in this reduced word for $w_yw$. As $w_x$ is the shortest element in $x$, to get $w_x$ we have to drop all letters in the reduced word for $w$ and some letters in the reduced word for $w_y$, so $w_x\le w_y$.


To prove (ii), recall that by a well-known property of the Bruhat order on $W$ we have either $\ell(s_\alpha w_x)=\ell(w_x)-1$ and $s_\alpha w_x\le w_x$, or $\ell(s_\alpha w_x)=\ell(w_x)+1$ and $s_\alpha w_x\ge w_x$, depending on whether $w_x^{-1}(\alpha)<0$ or $w_x^{-1}(\alpha)>0$. In the first case $s_\alpha w_x$ is the shortest element in $s_\alpha x$, as $\ell(s_\alpha w_x)=\ell_S(x)-1$ and we obviously always have $|\ell_S(s_\alpha x)-\ell_S(x)|\le1$ by definition of~$\ell_S$. Hence $w_{s_\alpha x}=s_\alpha w_x$. In the second case we have $s_\alpha x\ge x$ by~(i) and hence
$\ell_S(s_\alpha x)\ge\ell_S(x)$. Therefore either $\ell_S(s_\alpha x)=\ell_S(x)+1$, in which case $s_\alpha w_x$ is the shortest element in $s_\alpha x$ and hence $w_{s_\alpha x}=s_\alpha w_x$, or $\ell_S(s_\alpha x)=\ell_S(x)$, in which case $s_\alpha x=x$ as $s_\alpha x\ge x$.
\ep

Note that part (i) implies in particular that if $w_0$ is the longest element in $W$ then $x_0=w_0W_S$ is the largest, hence the longest, element in $X_S$. Part (ii) implies that if $x\in X_S$ and $w=s_{i_1}\dots s_{i_n}\in x$ is written in reduced form, then a reduced word for $w_x$ can be obtained from $s_{i_1}\dots s_{i_n}$ by dropping all letters $s_{i_j}$ such that $s_{i_n}\dots s_{i_{j+1}}s_{i_j}s_{i_{j+1}}\dots s_{i_n}\in W_S$.

\smallskip

Denote by $P_{++}(S^c)\subset P_+(S^c)$ the subset consisting of weights $\lambda$ such that $\lambda(H_\alpha)>0$ for every $\alpha\in S^c$. The following result is well-known in the case $S=\emptyset$, see \cite[Theorem 2.9]{BGG}.

\begin{proposition} \label{pBGG}
Let $\lambda\in P_{++}(S^c)$ and $x,y\in X_S$. Then $x\le y$ if and only if $V_\lambda(w_x\lambda)\subset (U_q\bb)V_\lambda(w_y\lambda)$.
\end{proposition}

\bp By virtue of Lemma \ref{lBruhat} the proof is essentially identical to that of \cite[Theorem 2.9]{BGG} for the case $S=\emptyset$, and proceeds along the following lines. Define a partial order on $X_S$ by declaring $x\preceq y$ iff $V_\lambda(w_x\lambda)\subset (U_q\bb)V_\lambda(w_y\lambda)$. Note that it is indeed a partial order, since the stabilizer of $\lambda$ in~$W$ is exactly $W_S$ by Chevalley's lemma, see~\cite[Proposition~2.72]{Knapp}. It is checked that this order has the properties
\begin{itemize}
\item[(i)] if $\ell_S(s_\alpha x)=\ell_S(x)+1$ for some $x\in X_S$ and $\alpha\in\Pi$ then $x\preceq s_\alpha x$;
\item[(ii)] if $x\preceq y$ and $\alpha\in\Pi$ then either $s_\alpha x\preceq y$ or $s_\alpha x\preceq s_\alpha y$.
\end{itemize}
It is proved then that the Bruhat order is the unique order on $X_S$ satisfying these properties.
\ep

As usual define an action of the Weyl group on $T$ by requiring $\lambda(w(t))=(w^{-1}\lambda)(t)$ for $\lambda\in P=\hat T$. For $z\in\T$ define an automorphism $\theta_z$ of $C(\bar\DD_q)$ by $\theta_z(Z_q)=\bar z Z_q$.

\begin{lemma} \label{ltconj}
For every simple root $\alpha\in\Pi$ and $t\in T$ we have $\pi_t\otimes\pi_{s_\alpha}=\theta_{\alpha(t)}\pi_{s_\alpha}\otimes\pi_{s_\alpha(t)}$ as homomorphisms $C(G_q)\to C(\bar\DD_{q^{d_\alpha}})$. In particular, for every $w\in W$ and $t\in T$ the kernels of $\pi_t\otimes\pi_w$ and $\pi_w\otimes\pi_{w^{-1}(t)}$ coincide.
\end{lemma}

\bp Consider first the case $G=SU(2)$. In this case the claim is that $\pi_t\otimes\rho_q=\theta_{t^2}\rho_q\otimes\pi_{t^{-1}}$ for $t\in T\cong\T$. This is immediate by definition~\eqref{erho} of $\rho_q$, as
$$
(\pi_t\otimes\iota)\Delta_q(\alpha)=t\alpha,\ \ (\pi_t\otimes\iota)\Delta_q(\gamma)=t^{-1}\gamma, \ \
(\iota\otimes\pi_{t^{-1}})\Delta_q(\alpha)=t^{-1}\alpha,\ \
(\iota\otimes\pi_{t^{-1}})\Delta_q(\gamma)=t^{-1}\gamma.
$$

Consider now the general case. Note that similarly to~\eqref{et} we have $(\pi_t\otimes\pi_w)(C^\lambda_{\zeta,\xi})=\pi_w(C^\lambda_{t^{-1}\zeta,\xi})$. Let $t
=\exp({2\pi i h})\in T$, $h\in i\mathfrak t\subset\h$. Write $h=c H_\alpha+h_2$ with $c\in\R$ and $h_2\in\ker\alpha$, and put $t_1=\exp({2\pi i c H_\alpha})$ and $t_2=\exp({2\pi i h_2})$. Then $t=t_1t_2$, $s_\alpha(t)=t_1^{-1}t_2$ and $\alpha(t)=\alpha(t_1)=e^{4\pi i c}$. The homomorphisms $\pi_{t_1}$ and $\pi_{s_\alpha}$ factor through $C(SU_{q^{d_\alpha}}(2))$, hence
$$
\pi_{t_1}\otimes\pi_{s_\alpha}=\theta_{\alpha(t_1)}\pi_{s_\alpha}\otimes\pi_{t_1^{-1}}=\theta_{\alpha(t)}\pi_{s_\alpha}\otimes\pi_{t_1^{-1}}.
$$
Observe next that since $t_2$ commutes with $E_\alpha, F_\alpha\in U_q\g$, the restrictions of the matrix coefficients~$C^\lambda_{\zeta, t_2\xi}$ and~$C^\lambda_{t_2^{-1}\zeta, \xi}$ to the algebra generated by $E_\alpha$, $F_\alpha$ and $K_\alpha$ coincide. In other words, $C^\lambda_{\zeta, t_2\xi}$ and~$C^\lambda_{t_2^{-1}\zeta, \xi}$ have the same images in $C(SU_{q^{d_\alpha}}(2))$. We then have
\begin{align*}
(\pi_t\otimes\pi_{s_\alpha})(C^\lambda_{\zeta,\xi})&=(\pi_{t_1}\otimes\pi_{s_\alpha})(C^\lambda_{t_2^{-1}\zeta,\xi})
=(\pi_{t_1}\otimes\pi_{s_\alpha})(C^\lambda_{\zeta,t_2\xi})=(\theta_{\alpha(t)}\pi_{s_\alpha}\otimes\pi_{t_1^{-1}})(C^\lambda_{\zeta,t_2\xi})\\
&=(\theta_{\alpha(t)}\pi_{s_\alpha}\otimes\pi_{t_1^{-1}t_2})(C^\lambda_{\zeta,\xi})
=(\theta_{\alpha(t)}\pi_{s_\alpha}\otimes\pi_{s_\alpha(t)})(C^\lambda_{\zeta,\xi}).
\end{align*}

If $w=s_{i_1}\dots s_{i_n}$ then by induction we get
$$
\pi_t\otimes\pi_w=(\theta_{z_1}\otimes\dots\otimes\theta_{z_n})\pi_w\otimes\pi_{w^{-1}(t)},
$$
where $z_k=(s_{i_1}\dots s_{i_{k-1}}\alpha_{i_k})(t)$. This gives the last statement in the formulation of the lemma.
\ep

For $q=1$ the above lemma implies that $t\Sigma_w=\Sigma_w w^{-1}(t)$. This slightly weaker statement can be deduced without any computations from the fact that every symplectic leaf intersects the normalizer of $T$ at a unique point. For $q<1$ the lemma implies that the representations $\pi_t\otimes\pi_w$ and $\pi_w\otimes\pi_{w^{-1}(t)}$ of $C(G_q)$ on $\ell^2(\Z_+)^{\otimes\ell(w)}$ are equivalent. This can also be easily proved by comparing the highest weights~\cite{So} of these representations.

\smallskip

For $z\in\T$ denote by $\chi_z$ the character of $C(\bar\DD_q)$ defined by $\chi_z(Z_q)=z$. Assume $\alpha\in\Pi$ and $c\in\R$. Put $t=\exp({2\pi i c H_\alpha})\in T$ and $z=e^{-2\pi i c}$. Then
\begin{equation} \label{echaracter}
\pi_t=\chi_z\pi_{s_\alpha}\ \ \hbox{on}\ \ C(G_q).
\end{equation}
Indeed, this is enough to check for $G=SU(2)$, and then this is immediate, since $\pi_t(\alpha_q)=\bar z$ and~$\pi_t(\gamma_q)=0$.

\begin{lemma} \label{lproductcells}
For every $1\le m\le m_0$ the ideal $J_m\subset C(G_q/K^{S,L}_q)$ is contained in the kernel of $\pi_{s_{i_1}}\otimes\dots\otimes\pi_{s_{i_n}}\otimes\pi_t$  for any $t\in T$ and any indices $i_1,\dots,i_n$ with $n\le m$.
\end{lemma}

\bp We will prove this in several steps. First note that $J_{m+1}\subset J_m$ for all $0\le m<m_0$. Indeed, let $w\in W^S$, $\ell(w)=m$, and $t\in T$. There exists $\alpha\in\Pi$ such that $s_\alpha w\in W^S$ and $\ell(s_\alpha w)=m+1$. (This follows e.g.~from Lemma~~\ref{lBruhat}(ii). Indeed, let $w_0$ be the longest element in~$W$. Write $w_0=s_{i_n}\dots s_{i_1}w$ with $n=\ell(w_0)-\ell(w)$. Then we can take $\alpha=\alpha_{i_k}$, where $k$ is the smallest number such that $w^{-1}s_{i_k}w\notin W_S$.) By \eqref{echaracter} we have $\pi_{w,t}=(\chi_1\otimes\iota)(\pi_{s_\alpha}\otimes\pi_{w,t})$. Therefore if $a\in J_{m+1}\subset\ker\pi_{s_\alpha w,t}=\ker(\pi_{s_\alpha}\otimes\pi_{w,t})$ then $a\in\ker \pi_{w,t}$. Thus $J_{m+1}\subset J_m$.

It follows that  if $a\in J_m$ for some $1\le m\le m_0$ then $\pi_{w,t}(a)=0$ for any $w\in W$ with $\ell(w)\le m$ and any $t\in T$. Indeed, write $w=w'w''$ with $w'\in W^S$ and $w''\in W_S$. Then, as already used in the proof of Theorem~\ref{treps}, we can assume that $\pi_w=\pi_{w'}\otimes\pi_{w''}$, so that $\pi_{w,t}(a)=\pi_{w',t}(a)\otimes 1^{\otimes\ell(w'')}$ by~\eqref{eDS}. Since $n=\ell(w')\le m$ and $J_m\subset J_n$, we have $\pi_{w',t}(a)=0$, hence also $\pi_{w,t}(a)=0$.

The next thing to observe is that if $a\in J_{m+1}$ for some $0\le m< m_0$ then $a\in\ker(\pi_{s_\alpha}\otimes\pi_{w,t})$
for any $t\in T$, $\alpha\in\Pi$ and $w\in W$ with $\ell(w)\le m$. If $\ell(s_\alpha w)=\ell(w)+1$, this is clearly true. Therefore assume that $\ell(s_\alpha w)=\ell(w)-1$. Put $w'=s_\alpha w$. Then $w=s_\alpha w'$ and we may assume that $\pi_w=\pi_{s_\alpha}\otimes\pi_{w'}$. Since $a\in J_{m+1}\subset J_m$, we have $a\in\ker(\pi_{s_\alpha}\otimes\pi_{w',\tau})$ for any $\tau\in T$.
Denote by $T_\alpha\subset T$ the set of elements $u$ of the form $u=\exp({2\pi i c H_\alpha})$, $c\in \R$.
By Lemma~\ref{ltconj} it follows that for any $u\in T_\alpha$ the element $a$ belongs to
$$
\ker(\pi_{s_\alpha,u}\otimes\pi_{w',t})=\ker(\pi_{s_\alpha}\otimes\pi_{w',w'^{-1}(u)t}).
$$
In other words, if $\varphi$ is a bounded linear functional on $C(G_q/K^{S,L}_q)$ of the form $\psi\pi_{w',t}$ then
$$
(\iota\otimes\varphi)\Delta_q(a)\in\ker \pi_{s_\alpha,u}\ \ \hbox{for any}\ \ u\in T_\alpha.
$$
Since the intersection of the kernels of the homomorphisms $\rho_q\otimes\pi_u\colon C(SU_q(2))\to C(\bar\DD_q)$, $u\in \T$, is zero, the intersection of the kernels of the homomorphisms $\pi_{s_\alpha,u}\colon C(G_q)\to C(\bar\DD_{q^{d_\alpha}})$, $u\in T_\alpha$, is exactly the kernel of the homomorphism $C(G_q)\to C(SU_{q^{d_\alpha}}(2))$. Therefore $(\iota\otimes\varphi)\Delta_q(a)$ is in the kernel of the latter homomorphism. Since $\pi_{s_\alpha}\otimes\pi_{s_\alpha}$ factors through the homomorphism $C(G_q)\to C(SU_{q^{d_\alpha}}(2))$, we conclude that
$$
(\iota\otimes\varphi)\Delta_q(a)\in\ker (\pi_{s_\alpha}\otimes\pi_{s_\alpha}).
$$
Since this is true for any $\varphi$ of the form $\psi\pi_{w',t}$, it follows that $a\in\ker(\pi_{s_\alpha}\otimes\pi_{s_\alpha}\otimes\pi_{w',t})$. As $\pi_{s_\alpha}\otimes\pi_{w',t}=\pi_{w,t}$, this proves the claim.

We now turn to the proof of the statement in the formulation. The proof is by induction on~$m$. For $m=1$ the result is already proved in the second paragraph. So assume the result is true for all numbers not bigger than $m<m_0$. Since $J_{m+1}\subset J_n$ for $n\le m$, it suffices to show that the kernel of $\pi_{s_{i_1}}\otimes\dots\otimes\pi_{s_{i_{m+1}}}\otimes\pi_t$ contains $J_{m+1}$ for any $i_1,\dots,i_{m+1}$ and $t\in T$. Let $a\in J_{m+1}$. Then by the previous paragraph $a\in\ker(\pi_{s_{i_1}}\otimes\pi_{w,t})$ for all $t\in T$ and all $w\in W$ with $\ell(w)\le m$ . Hence, for any bounded linear functional $\varphi$ on $C(G_q)$ of the form $\psi\pi_{s_{i_1}}$ we have $(\varphi\otimes\iota)\Delta_q(a)\in\ker\pi_{w,t}$. It follows that $(\varphi\otimes\iota)\Delta_q(a)\in J_m$. Hence, by the inductive assumption, $(\varphi\otimes\iota)\Delta_q(a)\in\ker(\pi_{s_{i_2}}\otimes\dots\otimes\pi_{s_{i_{m+1}}}\otimes\pi_t)$. Since this is true for any $\varphi$ of the form $\psi\pi_{s_{i_1}}$, we conclude that $a\in \ker(\pi_{s_{i_1}}\otimes\dots\otimes\pi_{s_{i_{m+1}}}\otimes\pi_t)$.
\ep

\bp[Proof of Theorem~\ref{tcompseries}] That $J_m\subset J_{m-1}$, follows from Lemma~\ref{lproductcells} (and was explicitly established in its proof). In particular, $J_{m_0}$ is contained in the kernel of every irreducible representation of~$C(G_q/K^{S,L}_q)$, hence $J_{m_0}=0$.

Let $1\le m\le m_0$. Consider the homomorphism
$$
\Theta_m\colon C(G_q/K^{S,L}_q)\to \bigoplus_{w\in W^S: \ell(w)=m}C(T/T_L;C(\bar\DD_{q_{i_1(w)}})\otimes\dots\otimes
C(\bar\DD_{q_{i_m(w)}}))
$$
defined by  $\Theta_m(a)_w(t)=\pi_{w,t}(a)$. The kernel of $\Theta_m$ is by definition the ideal $J_m$. Let $a\in J_{m-1}$, $t\in T$ and $w\in W^S$, $\ell(w)=m$. Let $w=s_{i_1}\dots s_{i_m}$ be the reduced expression used to define $\pi_w$. Let $1\le k\le m$ and $z=e^{-2\pi i c}\in\T$. Put $u=\exp({2\pi i cH_{i_k}})\in T$. Then applying $\chi_z\colon C(\bar \DD_{q_{i_k}})\to\C$ to the~$k$th factor of the image of $\pi_{w,t}$, by \eqref{echaracter} and Lemma~\ref{ltconj} we get
\begin{align*}
\ker((\iota\otimes\dots\otimes\chi_z\otimes\dots\otimes\iota)\pi_{w,t})
&=\ker(\pi_{s_{i_1}}\otimes\dots\otimes \pi_u\otimes\dots\otimes\pi_{s_{i_m}})\\
&=\ker(\pi_{s_{i_1}}\otimes\dots\otimes\pi_{s_{i_{k-1}}}\otimes\pi_{s_{i_{k+1}}}\otimes\dots\otimes\pi_{s_{i_m}}\otimes \pi_{u'}),
\end{align*}
where $u'=(s_{i_m}\dots s_{i_{k+1}})(u)$. Since $a\in J_{m-1}$, by Lemma~\ref{lproductcells} we thus see that $a$ is contained in the kernel of $(\iota\otimes\dots\otimes\chi_z\otimes\dots\otimes\iota)\pi_{w,t}$. Since this is true for all $z\in \T$, it follows that $\pi_{w,t}(a)$ is contained in the kernel of $\iota\otimes\dots\otimes\beta\otimes\dots\otimes\iota$, where $\beta\colon C(\bar\DD_{q_{i_k}})\to C(\T)$ is the homomorphism that maps $Z_{q_{i_k}}$ to the standard generator of $C(\T)$. The kernel of $\beta$ is by definition the ideal $C_0(\DD_{q_{i_k}})$. Therefore
$$
\pi_{w,t}(a)\in C(\bar\DD_{q_{i_1}})\otimes\dots\otimes C_0(\DD_{q_{i_k}})\otimes\dots
\otimes C(\bar\DD_{q_{i_m}}).
$$
Since this is true for every $k$, we conclude that $\pi_{w,t}(a)\in C_0(\DD_{q_{i_1}})\otimes\dots\otimes C_0(\DD_{q_{i_m}})$. Thus the image of $J_{m-1}$ under $\Theta_m$ is contained in
$$
\bigoplus_{w\in W^S: \ell(w)=m}C(T/T_L;C_0(\DD_{q_{i_1(w)}})\otimes\dots\otimes
C_0(\DD_{q_{i_m(w)}})).
$$
To see that this algebra is the whole image, we will consider separately the cases $q=1$ and $q<1$.

Assume $q=1$. In this case $J_{m-1}$ is the ideal of continuous functions on $G/K^{S,L}$ that vanish on the symplectic leaves of dimension $2m-2$. By the Stone-Weierstrass theorem it is then enough to show that for any two distinct points on the union of the leaves of dimension $2m$, there is a continuous function which vanishes on all leaves of dimension $2m-2$ and takes different nonzero values at these points. For this it suffices to know that the union of the leaves of dimension $\le 2m-2$ is a closed subset of $G/K^{S,L}$. This, in turn, is enough to check for $G/K^S$, which is the quotient of $G/K^{S,L}$ by an action of the compact group $T/T_L$. The result is well-known for $S=\emptyset$, see e.g.~Theorem~23 on p.~127 in~\cite{St}. Since the union of the symplectic leaves  of $G/K^S$ of dimension $\le 2m-2$ is the image of the union of the symplectic leaves  of $G/T$ of dimension $\le 2m-2$, we conclude that this set is closed for any $S$.

Turning to the case $q<1$, first we prove that $\pi_{w,t}(J_{m-1})\ne0$ for any $t\in T$ and $w\in W^S$ with $\ell(w)=m$. For this take any $\lambda\in P_{++}(S^c)$. Since $w$ cannot be smaller than any $v\in W^S$ with $\ell(v)=m-1$, by Proposition~\ref{pBGG} we see that $V_\lambda(w\lambda)$ is orthogonal to $(U_q\bb)V_\lambda(v\lambda)$, hence $\pi_{v,\tau}(C^\lambda_{w\lambda,\lambda})=\lambda(\tau)\pi_{v}(C^\lambda_{w\lambda,\lambda})=0$ for any $\tau\in T$ by Lemma~\ref{lsoib}(ii). Therefore $(C^\lambda_{w\lambda,\lambda})^*C^\lambda_{w\lambda,\lambda}\in J_{m-1}$. By Lemma~\ref{lsoib}(i) we also have $\pi_{w,t}((C^\lambda_{w\lambda,\lambda})^*C^\lambda_{w\lambda,\lambda})\ne0$.

Since $J_{m-1}$ is an ideal, it follows that the representations $\pi_{w,t}$ of $J_{m-1}$, with $w\in W^S$, $\ell(w)=m$ and $t\in T/T_L$, are irreducible and mutually inequivalent. In other words, the subalgebra $\Theta_m(J_{m-1})$ of the algebra
$$
\bigoplus_{w\in W^S: \ell(w)=m}C(T/T_L;C_0(\DD_{q_{i_1(w)}})\otimes\dots\otimes
C_0(\DD_{q_{i_m(w)}}))=\bigoplus_{w\in W^S: \ell(w)=m}C(T/T_L;K(\ell^2(\Z_+^m)))
$$
has the property that its projection to different fibers gives mutually inequivalent irreducible representations of~$\Theta_m(J_{m-1})$ on~$\ell^2(\Z_+^m)$. By the Stone-Weierstrass theorem for type I C$^*$-algebras we conclude that $\Theta_m(J_{m-1})$ coincides with the whole algebra.
\ep

\bigskip


\section{Topology on the spectrum} \label{s4}

In this section we will describe the Jacobson, or hull-kernel, topology on the spectrum of the type~I C$^*$-algebra $C(G_q/K^{S,L}_q)$ for $q\in(0,1)$. By Theorem~\ref{treps}, as a set the spectrum can be identified with $W^S\times T/T_L$.

To formulate the result it is convenient to use the description of the Bruhat order given in~\cite{BGG}. For $\sigma,w\in W$ and $\gamma\in\Delta_+$ we write $\sigma\xrightarrow{\gamma}w$ if $w=\sigma s_{\gamma}$ and $\ell(w)=\ell(\sigma)+1$ (note that in the notation of~\cite{BGG} this corresponds to $\sigma^{-1}\xrightarrow{\gamma}w^{-1}$). Then, for any $\sigma,w\in W$, we have $\sigma\le w$ if and only if there exist $\sigma_1,\dots,\sigma_k\in W$ and $\gamma_1,\dots,\gamma_k\in\Delta_+$ such that $\sigma\xrightarrow{\gamma_1}\sigma_1\xrightarrow{\gamma_2}\dots\xrightarrow{\gamma_k}\sigma_k=w$.

Assume $\sigma\le w$. For every path $\sigma\xrightarrow{\gamma_1}\sigma_1\xrightarrow{\gamma_2}\dots\xrightarrow{\gamma_k}\sigma_k=w$ consider the closed connected subgroup of~$T$ consisting of the elements of the form $\exp({i h_\beta})$ with $\beta\in\operatorname{span}_\R\{\gamma_1,\dots,\gamma_k\}$. Denote by~$T_{\sigma,w}$ the union of such groups for all possible paths. The closed sets~$T_{\sigma,w}$ clearly have the following multiplicative property, which will play an important role: if $\sigma\le v\le w$ then
\begin{equation}\label{eTmult}
T_{\sigma,v} T_{v,w}:=\{\tau t\mid \tau\in T_{\sigma,v},\ t\in T_{v,w}\}\subset T_{\sigma,w}.
\end{equation}

Recall that we denote by $\pi\colon G\to G/K^{S,L}$ the quotient map.

\begin{theorem} \label{ttopology}
Let $w\in W^S$ and $\Omega\subset T$. Then
\enu{i} the closure of the union of the symplectic leaves $\pi(\Sigma_w\tau)\subset G/K^{S,L}$, $\tau\in\Omega$, is the union of the leaves $\pi(\Sigma_\sigma t)$ such that $\sigma\in W^S$, $\sigma\le w$ and $t\in\bar\Omega T_{\sigma,w}T_L$;
\enu{ii} if $q\in(0,1)$, for any $\sigma\in W^S$ and $t\in T$, the kernel of the representation $\pi_{\sigma,t}$ contains the intersection of the kernels of the representations $\pi_{w,\tau}$  of $C(G_q/K^{S,L}_q)$, $\tau\in\Omega$, if and only if $\sigma\le w$ and $t\in\bar\Omega T_{\sigma,w}T_L$.
\end{theorem}

Therefore if for $q\in(0,1)$ we identify the spectrum of $C(G_q/K_q^{S,L})$ with the quotient of $G/K^{S,L}$ by the partition defined by its symplectic leaves (or in other words, with the quotient of $G/K^{S,L}$ by the right dressing action), then the Jacobson topology on the spectrum is exactly the quotient topology.

\smallskip

The proof is based on the following refinement of Lemma~\ref{lproductcells}. Recall that in the proof of that lemma we denoted by $T_\alpha\subset T$ the subgroup consisting of elements of the form $\exp({2\pi i c H_\alpha})$, $c\in \R$. We write $T_i$ for $T_{\alpha_i}$.

\begin{lemma} \label{lproductcells2}
Let $1\le i_1,\dots,i_n\le r$ and $t\in T$. Assume $a\in C(G_q)$ is such that
$\pi_{w,\tau}(a)=0$ for all $1\le k\le n$, $w=s_{i_{j_1}}\dots s_{i_{j_k}}$ with $1\le j_1<\dots <j_k\le n$ and $\tau\in T$ such that $\tau t^{-1}$ lies in the group generated by $(s_{i_{j_k}}s_{i_{j_{k-1}}}\dots s_{i_{j_l}})(T_{i_m})$ with $1\le l\le k$ and $j_{l-1}<m<j_l$ (we let $j_0=0$).
Then $a\in\ker(\pi_{s_{i_1}}\otimes\dots\otimes \pi_{s_{i_n}}\otimes\pi_t)$.
\end{lemma}

\bp The proof is by induction on $n$. For $n=1$ the statement is tautological. So assume $n>1$ and that the result is true for all numbers $<n$. By induction, exactly as in the proof of Lemma~\ref{lproductcells}, it suffices to show that
$(\pi_{s_{i_1}}\otimes\pi_{w,\tau})(a)=0$ for all $1\le k\le n-1$, $w=s_{i_{j_1}}\dots s_{i_{j_k}}$ with $2\le j_1<\dots <j_k\le n$ and $\tau\in T$ such that $\tau t^{-1}$ lies in the group generated by $(s_{i_{j_k}}s_{i_{j_{k-1}}}\dots s_{i_{j_l}})(T_{i_m})$ with $1\le l\le k$ and $j_{l-1}<m<j_l$ (with $j_0=1$).
To see that this is true, fix $k$, $w=s_{i_{j_1}}\dots s_{i_{j_k}}$ and $\tau$. If $\ell(s_{i_1} w)=\ell(w)+1$ then we may assume that $\pi_{s_{i_1}}\otimes\pi_w=\pi_{s_{i_1}w}$, and then the claim is part of the assumption. If $\ell(s_{i_1} w)=\ell(w)-1$, the claim is proved by the same argument as in the third paragraph of the proof of Lemma~\ref{lproductcells}, using that $\pi_{w,\tau w^{-1}(u)}(a)=0$ for any $u\in T_{i_1}$, which is true by assumption.
\ep

Our goal is to relate the groups in the above lemma to the sets $T_{\sigma,w}$.

\begin{lemma} \label{lpaths}
Assume $\sigma,w\in W$ and $\alpha\in\Pi$ are such that $\sigma<s_\alpha\sigma$ and $w<s_\alpha w$. Then for any path $\sigma\xrightarrow{\gamma_1}\dots\xrightarrow{\gamma_k}w$ there exists a path $s_\alpha\sigma\xrightarrow{\gamma_1'}\dots\xrightarrow{\gamma_k'}s_\alpha w$ such that the group generated by $\gamma_1',\dots,\gamma_k'$ coincides with the group generated by $\gamma_1,\dots,\gamma_k$.
\end{lemma}

\bp The proof is by induction on $k=\ell(w)-\ell(\sigma)$. Put $v=\sigma\gamma_1$. Consider the two possible cases.

Assume first that $v<s_\alpha v$. By the inductive assumption there exists a path $s_\alpha v\xrightarrow{\gamma_2'}\dots\xrightarrow{\gamma_k'}s_\alpha w$ such that the group generated by $\gamma_2',\dots,\gamma_k'$ coincides with the group generated by $\gamma_2,\dots,\gamma_k$. Then $s_\alpha\sigma\xrightarrow{\gamma_1}s_\alpha v\xrightarrow{\gamma_2'}\dots\xrightarrow{\gamma_k'}s_\alpha w$ is the required path.

Assume now that $s_\alpha v<v$. Then $v=s_\alpha \sigma$. In particular, $\gamma_1=\sigma^{-1}(\alpha)$. Consider the path $s_\alpha\sigma\xrightarrow{\gamma_2}\dots\xrightarrow{\gamma_k}w \xrightarrow{w^{-1}(\alpha)} s_\alpha w$. Since $w=s_\alpha\sigma s_{\gamma_2}\dots s_{\gamma_k}$, we have $w^{-1}(\alpha)=-(s_{\gamma_k}\dots s_{\gamma_2}\sigma^{-1})(\alpha)$. Therefore the groups generated by $\gamma_1=\sigma^{-1}(\alpha),\gamma_2,\dots,\gamma_k$ and by $\gamma_2,\dots,\gamma_k, w^{-1}(\alpha)$ coincide.
\ep

Note that the proof actually shows that as the sequence $\gamma_1',\dots,\gamma_k'$ we can take $$\gamma_1,\dots,\gamma_{i-1},\gamma_{i+1},\dots,\gamma_k,w^{-1}(\alpha),$$ where $i$ is the first number such that $s_\alpha\sigma s_{\gamma_1}\dots s_{\gamma_{i-1}}=\sigma s_{\gamma_1}\dots s_{\gamma_i}$ (if there is no such number then the sequence is $\gamma_1,\dots,\gamma_k$).

\begin{lemma} \label{lTsets}
Let $w=s_{i_1}\dots s_{i_n}$ be written in reduced form, and consider $\sigma\le w$.
\enu{i} Assume $\sigma =s_{i_{j_1}}\dots s_{i_{j_k}}$ for some $0\le k\le n$ and $1\le j_1<\dots< j_k\le n$. Let $\Gamma\subset P$ be the group generated by $(s_{i_{j_k}}s_{i_{j_{k-1}}}\dots s_{i_{j_l}})(\alpha_{i_m})$ with $1\le l\le k+1$ and $j_{l-1}<m<j_l$ (we let $j_0=0$ and $j_{k+1}=n+1$). Then $\Gamma$ coincides with the group generated by the elements $(s_{i_n}\dots s_{i_{m+1}})(\alpha_{i_m})$ such that $m\notin\{j_1,\dots,j_k\}$.
\enu{ii} Under the assumptions of {\rm (i)}, there exists a path $\sigma\xrightarrow{\gamma_1}\dots\xrightarrow{\gamma_p}w$ such that $\Gamma$ is contained in the group generated by $\gamma_1,\dots,\gamma_p$, and these two groups coincide if the expression $\sigma=s_{i_{j_1}}\dots s_{i_{j_k}}$ is reduced.
\enu{iii} For any path $\sigma\xrightarrow{\gamma_1}\dots\xrightarrow{\gamma_{n-k}}w$ there exist $1\le j_1<\dots< j_k\le n$ such that $\sigma =s_{i_{j_1}}\dots s_{i_{j_k}}$ and the group $\Gamma$ defined as in {\rm (i)} coincides with the group generated by $\gamma_1,\dots,\gamma_{n-k}$.
\end{lemma}

\bp (i) The proof is by induction on $n$. For $n=1$ the statement is tautological. Assume $n>1$. If $j_k=n$ then the result is immediate by the inductive assumption. Assume $j_k<n$. Let $w'=s_{i_1}\dots s_{i_{n-1}}$ and $\Gamma'$ be the group defined similarly to $\Gamma$ by the elements $\sigma\le w'$. Then $\Gamma$ is generated by $\Gamma'$ and~$\alpha_{i_n}$, hence by $s_{{i_n}}(\Gamma')$ and $\alpha_{i_n}$. Since by the inductive assumption $\Gamma'$ is generated by the elements $(s_{i_{n-1}}\dots s_{i_{m+1}})(\alpha_{i_m})$ such that $m\le n-1$ and $m\notin\{j_1,\dots,j_k\}$, we get the result.

\smallskip

(ii) The proof is again by induction on $n$. For $n=1$ the statement is trivial. So assume $n>1$. We may assume that the word $s_{i_{j_1}}\dots s_{i_{j_k}}$ is reduced. Indeed, if it is not, then a reduced expression for~$\sigma$ can be obtained by dropping some letters in the word $s_{i_{j_1}}\dots s_{i_{j_k}}$, see  Lemma~21(c) in Appendix to~\cite{St}. By the description of $\Gamma$ given in (i) this can only increase the group~$\Gamma$. Consider two cases.

Assume $k\ge1$ and $j_1=1$. Put $\sigma'=s_{i_{j_2}}\dots s_{i_{j_k}}$ and $w'=s_{i_2}\dots s_{i_n}$. By the inductive assumption there exists a path $\sigma'\xrightarrow{\gamma_1'}\dots\xrightarrow{\gamma_{n-k}'}w'$ such that $\Gamma$ coincides with the group generated by $\gamma_1',\dots,\gamma_{n-k}'$. By Lemma~\ref{lpaths} we can then find a path $\sigma\xrightarrow{\gamma_1}\dots\xrightarrow{\gamma_{n-k}}w$ such that the group generated by $\gamma_1,\dots,\gamma_{n-k}$ coincides with the group generated by $\gamma_1',\dots,\gamma_{n-k}'$.

Assume now that either $k=0$ or $j_1>1$. Let $v=s_{i_2}\dots s_{i_n}$. Then by (i) and the inductive assumption there exists a path $\sigma\xrightarrow{\gamma_1}\dots\xrightarrow{\gamma_{n-k-1}}v$ such that $\Gamma$ coincides with the group generated by $\gamma_1,\dots,\gamma_{n-k-1}$
and $v^{-1}(\alpha_{i_1})$. Therefore we can take $\gamma_{n-k}=v^{-1}(\alpha_{i_1})$, so that we get a path $\sigma\xrightarrow{\gamma_1}\dots\xrightarrow{\gamma_{n-k-1}}v\xrightarrow{\gamma_{n-k}}w$.

\smallskip

(iii) By \cite[Proposition~2.8(c)]{BGG}, given a path $\sigma\xrightarrow{\gamma_1}\dots\xrightarrow{\gamma_{n-k}}w$ there exist uniquely defined numbers $p_1,\dots,p_{n-k}$ such that $\sigma\gamma_1\dots\gamma_{n-k-l}$ is obtained from $s_{i_1}\dots s_{i_n}$ by dropping the letters $s_{i_{p_1}},\dots,s_{i_{p_l}}$. Let $\{j_1<\dots <j_k\}$ be the complement  of $\{p_1,\dots,p_{n-k}\}$ in $\{1,\dots,n\}$. It remains to show that the group~$\Gamma$ is generated by $\gamma_1,\dots,\gamma_{n-k}$. Once again the proof is by induction on~$n$. Put $p=p_1$. Consider the element $w'=\sigma\gamma_1\dots\gamma_{n-k-1}=s_{i_1}\dots\hat s_{i_{p}}\dots s_{i_n}$. Let $\Gamma'$ be the group defined similarly to $\Gamma$ by the elements $\sigma\le w'$. By the inductive assumption it is generated by $\gamma_1,\dots,\gamma_{n-k-1}$. We also have $\gamma_{n-k}=(s_{i_n}\dots s_{i_{p+1}})(\alpha_{i_p})$. By part (i) the group $\Gamma'$ is generated by the elements $(s_{i_n}\dots s_{i_{m+1}})(\alpha_{i_m})$ such that  $m\in\{p_2,\dots,p_{n-k}\}$ and $m>p$ and the elements $(s_{i_n}\dots\hat s_{i_p}\dots s_{i_{m+1}})(\alpha_{i_m})$ such that $m\in\{p_2,\dots,p_{n-k}\}$ and $m<p$. Since $\gamma_{n-k}=(s_{i_n}\dots s_{i_{p+1}})(\alpha_{i_p})$, for $m<p$ we have
$$
s_{i_n}\dots\hat s_{i_p}\dots s_{i_{m+1}}=s_{\gamma_{n-k}}s_{i_n}\dots s_{i_{m+1}}.
$$
Therefore the group generated by $\gamma_{n-k}$ and $\Gamma'$ coincides with the group generated by the elements $(s_{i_n}\dots s_{i_{m+1}})(\alpha_{i_m})$ such that  $m\in\{p,p_2,\dots,p_{n-k}\}$, which is exactly the group $\Gamma$.
\ep

The previous lemma shows that the collection $X_{\sigma,w}$ of groups generated by $\gamma_1,\dots,\gamma_{n-k}$, where $n=\ell(w)$ and $k=\ell(\sigma)$, for all possible paths $\sigma\xrightarrow{\gamma_1}\dots\xrightarrow{\gamma_{n-k}}w$, can be described as follows. Fix a reduced decomposition $w=s_{i_1}\dots s_{i_n}$. For every sequence $1\le j_1<\dots< j_k\le n$ such that $\sigma=s_{j_1}\dots s_{j_k}$ consider the group generated by the elements $(s_{i_n}\dots s_{i_{m+1}})(\alpha_{i_m})$ such that $m\notin\{j_1,\dots,j_k\}$. Then $X_{\sigma,w}$ consists of such groups for all possible $j_1,\dots, j_k$.

In the particular case $\sigma=e$ this implies that $X_{e,w}$ consists of just one group, and for any reduced decomposition $w=s_{i_1}\dots s_{i_n}$ this group is generated by the elements $(s_{i_n}\dots s_{i_{m+1}})(\alpha_{i_m})$, $1\le m\le n$ (or equivalently, by the elements $\alpha_{i_1},\dots,\alpha_{i_n}$). That the latter group is independent of the reduced decomposition is well-known. In fact, the set of elements $(s_{i_n}\dots s_{i_{m+1}})(\alpha_{i_m})$, $1\le m\le n$, is exactly $\Delta_+\cap w^{-1}\Delta_-$, see e.g.~Corollary 2 to Proposition VI.6.17 in~\cite{Bour}. Therefore the set $T_{e,w}$ is the group consisting of the elements $\exp(ih_\beta)$ with $\beta\in\operatorname{span}_\R(\Delta_+\cap w^{-1}\Delta_-)$. It would be interesting to have a geometric description of $X_{\sigma,w}$ and $T_{\sigma,w}$ for all $\sigma\le w$.

\smallskip

The following lemma improves the main part of the proof of Theorem~\ref{tcompseries}.

\begin{lemma}
Let $t\in T$ and let $w=s_{i_1}\dots s_{i_n}\in W^S$ be written in the reduced form used to define~$\pi_w$. Assume $a\in C(G_q/K^{S,L}_q)$ is such that $\pi_{\sigma,\tau}(a)=0$ for all $\sigma\in W^S$ such that $\sigma<w$ and all $\tau\in tT_{\sigma,w}$. Then
$$
\pi_{w,t}(a)\in C_0(\DD_{q_{i_1}})\otimes\dots\otimes
C_0(\DD_{q_{i_{n}}}).
$$
\end{lemma}

\bp As in the proof of Theorem~\ref{tcompseries} it suffices to check that for any $z\in\T$ and $1\le m\le n$, applying $\chi_z$ to the $m$th factor of
$\pi_{w,t}(a)\in C(\bar\DD_{q_{i_1}})\otimes\dots\otimes
C(\bar\DD_{q_{i_{n}}})$ we get zero. Assume $z=e^{-2\pi ic}$ and put $u=\exp({2\pi i cH_{i_m}})$. Then exactly as in the proof of Theorem~\ref{tcompseries} we have to check that
$$
a\in\ker(\pi_{s_{i_1}}\otimes\dots\otimes\pi_{s_{i_{m-1}}}\otimes\pi_{s_{i_{m+1}}}\otimes\dots\otimes\pi_{s_{i_n}}\otimes \pi_{tu'}),
$$
where $u'=(s_{i_n}\dots s_{i_{m+1}})(u)$. For this, by Lemmas~\ref{lproductcells2} and~\ref{lTsets}, it suffices to check that $\pi_{\sigma,\tau}(a)=0$ for every $\sigma=s_{i_{j_1}}\dots s_{i_{j_k}}$ such that $m\notin\{j_1,\dots,j_k\}$ and all $\tau\in tT_{\sigma,w}$. If $\sigma\in W^S$, this is true by assumption. Otherwise write $\sigma=\sigma'\sigma''$ with $\sigma'\in W^S$ and $\sigma''\in W_S$. Then we may assume that $\pi_\sigma=\pi_{\sigma'}\otimes\pi_{\sigma''}$ and then by~\eqref{eDS} we have $\pi_{\sigma,\tau}(a)=\pi_{\sigma',\tau}(a)\otimes 1^{\otimes\ell(\sigma'')}$. Since $\sigma'\le\sigma<w$, we have $T_{\sigma,w}\subset T_{\sigma',w}$, and hence $\pi_{\sigma',\tau}(a)=0$ by assumption. Therefore we still get $\pi_{\sigma,\tau}(a)=0$.
\ep

\bp[Proof of Theorem~\ref{ttopology}] The main part of the argument works for all $q\in(0,1]$. Namely, we will show that the kernel of $\pi_{\sigma,t}$ contains the intersection of the kernels of the representations $\pi_{w,\tau}$ of~$C(G_q/K^{S,L}_q)$, $\tau\in\Omega$, if and only if $\sigma\le w$ and $t\in\bar\Omega T_{\sigma,w}T_L$.

\smallskip

Assume $\sigma\nleq w$ and $t\in T$. For $q=1$ it is known that the set $\cup_{v\le w}\Sigma_{v}T$ is closed in $G$, see again \cite[Theorem~23]{St} or \cite[Theorem~2.11]{BGG}, hence $\cup_{W^S\ni v\le w}\pi(\Sigma_{v}T)$ is closed in $G/K^{S,L}$ and does not intersect $\pi(\Sigma_\sigma t)$. For $q<1$, using Proposition~\ref{pBGG} and Lemma~\ref{lsoib}, for any $\lambda\in P_{++}(S^c)$ we get $\pi_{\sigma,t}((C^\lambda_{\sigma\lambda,\lambda})^*C^\lambda_{\sigma\lambda,\lambda})\ne0$ and $\pi_{w,\tau}((C^\lambda_{\sigma\lambda,\lambda})^*C^\lambda_{\sigma\lambda,\lambda})=0$ for all $\tau\in T$.

\smallskip

Assume $\sigma\le w$. Let $w=s_{i_1}\dots s_{i_n}$ be the reduced expression used to define $\pi_w$. For every $\tau\in T_{\sigma,w}$, by Lemma~\ref{lTsets} we can find $j_1<\dots <j_k$ such that $\sigma=s_{i_{j_1}}\dots s_{i_{j_k}}$ is reduced and $\tau$ has the form $\exp({i h_\beta})$ with $\beta$ lying in the real span of the elements $(s_{i_n}\dots s_{i_{m+1}})(\alpha_{i_m})$, $m\notin\{j_1,\dots,j_k\}$. Using Lemma~\ref{ltconj} we conclude that by applying the characters $\chi_{z_m}$ to the $m$th factor of $\pi_{w,\tau'}$ for appropriate numbers $z_{i_m}\in\T$ we can factor $\pi_{\sigma,\tau\tau'}$ through $\pi_{w,\tau'}$ for any $\tau'\in T$. In other words, if $t\in \Omega T_{\sigma,w}T_L$ then the kernel of $\pi_{\sigma,t}$ contains the kernel of $\pi_{w,\tau}$ for some $\tau\in\Omega$. In particular, the intersection of the kernels of $\pi_{w,\tau}$, $\tau\in\Omega$, is contained in the kernel of $\pi_{\sigma,t}$ for any $t\in\Omega T_{\sigma,w}T_L$. Since clearly the map $T\ni t\mapsto\pi_{\sigma,t}(a)$ is continuous for every $a\in C(G_q/K^{S,L}_q)$, the same is true for $t\in\bar\Omega T_{\sigma,w}T_L$.

\smallskip

Finally, assume that $\sigma\le w$, but $t\notin\bar\Omega T_{\sigma,w}T_L$. Let $m=\ell(\sigma)$ and $n=\ell(w)$. By Theorem~\ref{tcompseries} we can find $a_m\in J_{m-1}$ such that $\pi_{\sigma,t}(a_m)\ne0$, $\pi_{\sigma,\tau}(a_m)=0$ for all $\tau \in\bar\Omega T_{\sigma,w}T_L$, and $\pi_{\sigma',\tau}(a_m)=0$ for all $\tau\in T$ and $\sigma'\in W^S$ such that $\sigma'\ne\sigma$, $\ell(\sigma')=m$. If $n=m$ this already shows that the intersection of the kernels of $\pi_{w,\tau}$, $\tau\in\Omega$, is not contained in the kernel of $\pi_{\sigma,t}$. So assume $n>m$. We will construct by induction elements $a_k\in C(G_q/K^{S,L}_q)$, $m+1\le k\le n$, such that $a_k-a_{k-1}\in J_{k-1}$ for $m+1\le k\le n$, and $\pi_{v,\tau}(a_k)=0$ for $m\le k\le n$ for any $v\in W^S$ such that $v\le w$, $\ell(v)\le k$, and any $\tau \in\bar\Omega T_{v,w}T_L$. Note that $a_m$ satisfies the last requirement for $k=m$, since $a_m\in J_{m-1}$ and hence $\pi_{v,\tau}(a_m)=0$ if $\ell(v)\le m-1$ for any $\tau\in T$.

Assume $a_k$ is constructed. Let $v\in W^S$, $v\le w$, $\ell(v)=k+1$. By construction we have $\pi_{\sigma',\tau}(a_k)=0$ for any $\sigma'\in W^S$ such that $\sigma'<v$ and any $\tau\in \bar\Omega T_{\sigma',w}T_L$. Since $T_{\sigma',v}T_{v,w}\subset T_{\sigma',w}$ by~\eqref{eTmult}, by Lemma~\ref{lproductcells2} it follows that
$$
\pi_{v,\tau}(a_k)\in C_0(\DD_{q_{i_1(v)}})\otimes\dots\otimes
C_0(\DD_{q_{i_{k+1}(v)}})
$$
for any $\tau\in \bar\Omega T_{v,w}T_L$, where $v=s_{i_1(v)}\dots s_{i_{k+1}(v)}$ is the reduced decomposition used to define $\pi_{v}$. By Theorem~\ref{tcompseries} we can find $b\in J_k$ such that $\pi_{v,\tau}(a_k)=\pi_{v,\tau}(b)$ for all $v\in W^S$, $v\le w$, $\ell(v)=k+1$, and all $\tau\in \bar\Omega T_{v,w}T_L$. We then take $a_{k+1}=a_k-b$.

By construction we have $\pi_{w,\tau}(a_n)=0$ for all $\tau\in \bar\Omega T_L$.  As $a_n-a_m\in J_m$, we also have $\pi_{\sigma,t}(a_n)=\pi_{\sigma,t}(a_m)\ne0$.

\smallskip

This finishes the proof of (ii). For $q=1$ what we have proved means that a leaf $\pi(\Sigma_\sigma t)$ is contained in the closure of the union of the leaves $\pi(\Sigma_w\tau)$, $\tau\in\Omega$, if and only if $\sigma\le w$ and $t\in\bar\Omega T_{\sigma,w}T_L$. To establish (i) it remains to note  that since the symplectic leaves are  orbits of the right dressing action, the closure of the union of the leaves $\pi(\Sigma_w\tau)$, $\tau\in\Omega$, consists of entire leaves, so if a leaf $\pi(\Sigma_\sigma t)$ is not contained in this closure, it does not intersect it.
\ep

\bigskip


\section{Strict deformation quantization} \label{s5}

In this section we will consider the family of C$^*$-algebras $C(G_q/K^{S,L}_q)$. To distinguish elements of different algebras we will use upper and lower indices $q$. Indices corresponding to $q=1$ will often be omitted.

In \cite{NT6} we showed that the family $(C(G_q))_q$ has a canonical structure of a continuous field of C$^*$-algebras. It is defined as follows.  For every $q\in(0,1]$ choose a $*$-isomorphism $\varphi^q\colon\U(G_q)\to\U(G)$ extending the canonical identifications of the centers. In other words, for every $\lambda\in P_+$ choose a $*$-isomorphism $\varphi^q_\lambda\colon B(V^q_\lambda)\to B(V_\lambda)$. Then upon identifying $\U(G_q)$ with $\prod_{\lambda\in P_+}B(V^q_\lambda)$ the isomorphism $\varphi^q$ is given by $(\varphi^q_\lambda)_\lambda$. The family of isomorphisms $\{\varphi^q\}_q$ is called continuous if the maps $q\mapsto \varphi^q(X^q)\in\U(G)=\C[G]^*$ are $\sigma(\U(G),\C[G])$-continuous for $X^q=E_i^q,F_i^q,h_i$; in other words, for every $\lambda\in P_+$ and $\xi\in V_\lambda$, the maps $q\mapsto \varphi^q(X^q)\xi\in V_\lambda$ are continuous. By \cite[Lemma~1.1]{NT6} there always exists a continuous family of $*$-isomorphisms such that $\varphi^1=\iota$. Fix such a family and consider the dual maps $\hat\varphi^q\colon\C[G]\to\C[G_q]$. They are coalgebra isomorphisms. Then by \cite[Proposition~1.2]{NT6} the family $(C(G_q))_q$ has a unique structure of a continuous field of C$^*$-algebras such that for every $a\in\C[G]$ the section $q\mapsto\hat\varphi^q(a)\in C(G_q)$ is continuous, and this structure does not depend on the choice of a continuous family of isomorphisms. The proof is based on two results, which we will now recall as they both play an important role in what follows.

The first one, exploited in one way or another in all cases where a continuous field has been constructed~\cite{Bau}, \cite{Bl}, \cite{Na}, \cite{Sh0}, is that in view of the classification of irreducible representations the key step is to prove continuity for quantum disks. We include a sketch of a by now standard proof for the reader's convenience.

\begin{lemma} \label{ldiskcont}
The family of C$^*$-algebras $C(\bar\DD_q)$, $q\in[0,1]$, has a unique structure of a continuous field of C$^*$-algebras such that $q\mapsto Z_q$ is a continuous section. The continuous field $(C(\bar\DD_q))_{q\in[0,1)}$ is isomorphic to the constant field with fiber $\TT$.
\end{lemma}

\bp Consider the universal unital C$^*$-algebra $A$ generated by two elements $Z$ and $Q$ such that
$$
1-Z^*Z=Q^2(1-ZZ^*),\ \ QZ=ZQ, \ \ \|Z\|\le1, \ \ 0\le Q\le1.
$$
For $q\in[0,1]$ let $I_q\subset A$ be the ideal generated by $Q-q1$. Put $A_q=A/I_q$ and denote by $\pi_q$ the quotient map $A\to A_q$. Since $Q$ is in the center of $A$, the function $q\mapsto\|\pi_q(a)\|$ is automatically upper semicontinuous for every $a\in A$. It is also clear that $A_q\cong C(\bar\DD_q)$. Therefore to prove the lemma we just have to check that the functions $q\mapsto\|\pi_q(a)\|$ are lower semicontinuous. For this define states $\psi_q$ on $A_q$ as follows. The state $\psi_1$ on $A_1\cong C(\bar\DD)$ is given by the normalized Lebesgue measure on the unit disk. For $q<1$ the state $\psi_q$ on $A_q\cong\TT$ is defined by
$$
\psi_q(a)=(1-q^2)\sum^\infty_{n=0}q^{2n}(ae_n,e_n)\ \ \hbox{for}\ \ a\in\TT.
$$
It is not difficult to check that the family $(\psi_q)_q$ is continuous in the sense that the map $q\mapsto\psi_q(\pi_q(a))$ is continuous for every $a\in A$. Since the states are faithful, the corresponding GNS-representations~$\pi_{\psi_q}$ are faithful as well, which implies that the functions $q\mapsto\|\pi_q(a)\|=\|\pi_{\psi_q}(\pi_q(a))\|$ are lower semicontinuous.

The last statement in the formulation is immediate from the explicit isomorphism $C(\bar\DD_q)\cong\TT$.
\ep

The second result is that under the homomorphism $C(G_q)\to C(SU_{q^{d_\alpha}}(2))$ corresponding to the simple root $\alpha$ the image of $\hat\varphi^q(a)$ is a polynomial in the standard generators of $\C[SU_{q^{d_\alpha}}(2)]$ with coefficients that are continuous in $q$. This is a consequence of the following lemma, which we formulate in a more general setting needed later. The group $\tilde K^S$ is simply connected, semisimple (if nontrivial) and compact, and its set of dominant integral weights can be identified with $P_+(S)$. So for the same reasons as for $G$ we have a continuous family of $*$-isomorphisms $\U(\tilde K^S_q)\to\U(\tilde K^S)$ extending the identification of the centers of these algebras with the algebra of functions on $P_+(S)$. Slightly more generally, as $T_L=T_{P(S^c)}\times (P(S)+L)^\perp$, from Proposition~\ref{ppoissonsub}(ii) we get $K^{S,L}=\tilde K^S\times(P(S)+L)^\perp$, and therefore the irreducible representations of $K^{S,L}$ are classified by $P_+(S)\times P/(P(S)+L)=P_+(S)\times P(S^c)/L$. It follows then that the irreducible corepresentations of $C(K^{S,L}_q)$ are classified by a subset of $P_+(S)\times P(S^c)/L$, and the compact quantum group $K^{S,L}_q$ is a quotient of $\tilde K^S_q\times(P(S)+L)^\perp$. Therefore there exist injective $*$-homomorphisms $\psi^q\colon\U(K^{S,L}_q)\to \U(K^{S,L})$ extending the embeddings of the equivalence classes of irreducible corepresentations of $C(K^{S,L}_q)$ into $P_+(S)\times P(S^c)/L$. Then we say that a family $\{\psi^q\}_q$ of such homomorphisms is continuous if the maps $q\mapsto \varphi^q(X^q)\in\U(K^{S,L})$ are $\sigma(\U(K^{S,L}),\C[K^{S,L}])$-continuous for $X^q=E_i^q,F_i^q$ with $\alpha_i\in S$, for $X^q=h_\beta$ with $\beta\in L^\perp$, and for $X^q=t\in T_L$.

\begin{lemma}
We have $K^{S,L}_q=\tilde K^S_q\times (P(S)+L)^\perp$ for all $q\in(0,1)$, hence there exists a continuous family $\{\psi^q\colon\U(K^{S,L}_q)\to\U(K^{S,L})\}_{q\in(0,1]}$ of $*$-isomorphisms with $\psi^1=\iota$. For any such a family~$\{\psi^q\}_q$, there exists a continuous family $\{\varphi^q\colon\U(G_q)\to\U(G)\}_{q\in(0,1]}$ of $*$-isomorphisms with $\varphi^1=\iota$ and such that $\varphi^q=\psi^q$ on $\U(K^{S,L}_q)$.
\end{lemma}

This is established in the course of the proof of \cite[Proposition~1.2]{NT6} in the particular case when $K^{S,L}=SU(2)_\alpha$ is the subgroup corresponding to a simple root $\alpha$, so $S=\{\alpha\}$ and $L=P(S^c)$. The general case is proved in the same way.

\begin{proposition}
The family of C$^*$-algebras $(C(G_q/K^{S,L}_q))_{q\in(0,1]}$ is a continuous subfield of the continuous field $(C(G_q))_{q\in(0,1]}$.
\end{proposition}

\bp We have to check that the family $(C(G_q/K^{S,L}_q))_{q\in(0,1]}$ has enough continuous sections. It suffices to show that for every $q_0\in(0,1]$ and $a\in\C[G_{q_0}/K^{S,L}_{q_0}]$ there exists a continuous section $q\mapsto a(q)$ of the field $(C(G_q))_{q\in(0,1]}$ such that $a(q)\in C(G_q/K^{S,L}_q)$ for all $q$ and $a(q_0)=a$. Let  $\{\varphi^q\colon\U(G_q)\to\U(G)\}_{q\in(0,1]}$ be a continuous family of $*$-isomorphisms with $\varphi^1=\iota$ such that $\varphi^q(\U(K^{S,L}_q))=\U(K^{S,L})$. Then, since $\C[G_q/K^{S,L}_q]$ consists of the elements $b\in\C[G_q]$ such that $(\omega\otimes\iota)\Delta_q(b)=\hat\eps_q(\omega)b$ for all $\omega\in\U(K^{S,L}_q)$, $\hat\varphi^q$ is a coalgebra isomorphism and $\hat\eps_q=\hat\eps\varphi^q$, we have $\hat\varphi^q(\C[G/K^{S,L}])=\C[G_q/K^{S,L}_q]$. Since the section $q\mapsto\hat\varphi^q(b)$ is continuous for any $b\in\C[G/K^{S,L}]$ by definition of the continuous field structure on $(C(G_q))_{q\in(0,1]}$, we get the result.
\ep

For $0<a<b\le1$ denote by $\Gamma((C(G_q/K^{S,L}_q))_{q\in[a,b]})$ the C$^*$-algebra of continuous sections of the field $(C(G_q/K^{S,L}_q))_{q\in[a,b]}$. Let $\{\varphi^q\colon\U(G_q)\to\U(G)\}_{q\in(0,1]}$ be a continuous family of $*$-isomorphisms. Denote by $\Gamma_{alg}((\C[G_q])_{q\in[a,b]})$ the space of sections of the form $q\mapsto\sum^{n}_{i=1}f_i(q)\hat\varphi^q(a_i)$, where $n\in\N$, $a_i\in\C[G]$ and $f_i\in C[a,b]$. By \cite[Remark~1.3]{NT6} the space $\Gamma_{alg}((C[G_q])_{q\in[a,b]})$ does not depend on the choice of $\varphi^q$ and forms a dense $*$-subalgebra of $\Gamma((C(G_q))_{q\in[a,b]})$. Put
$$
\Gamma_{alg}((\C[G_q/K^{S,L}_q])_{q\in[a,b]})=\Gamma((C(G_q/K^{S,L}_q))_{q\in[a,b]})\cap\Gamma_{alg}((\C[G_q])_{q\in[a,b]}).
$$
Then $\Gamma_{alg}((\C[G_q/K^{S,L}_q])_{q\in[a,b]})$ is a dense involutive $C[a,b]$-subalgebra of $\Gamma((C(G_q/K^{S,L}_q))_{q\in[a,b]})$.

\smallskip

Recall that $\C[G/K^{S,L}]$ is a Poisson algebra with Poisson bracket defined by \eqref{epoissonbracket}.

\begin{theorem} \label{tquantization}
Assume $\eps\in(0,1)$. Then for any $a,b\in \Gamma_{alg}((\C[G_q/K^{S,L}_q])_{q\in[\eps,1]})$ we have
$$
\lim_{h\downarrow0}\frac{[a(e^{-h}),b(e^{-h})]}{ih}=\{a(1),b(1)\}.
$$
\end{theorem}

Here by the limit we mean that for some (equivalently, for any) $c\in \Gamma((C(G_q/K^{S,L}_q))_{q\in[\eps,1]})$ such that $\{a(1),b(1)\}=c(1)$ we have
$$
\lim_{h\downarrow0}\|[a(e^{-h}),b(e^{-h})]/{ih}-c(e^{-h})\|=0.
$$
Another way of formulating this theorem is to say that the section $c$ defined by $c(1)=\{a(1),b(1)\}$ and $c(e^{-h})=[a(e^{-h}),b(e^{-h})]/{ih}$ for $h>0$, is continuous.

\bp[Proof of Theorem~\ref{tquantization}] We may assume that $K^{S,L}$ is trivial. Let $\{\varphi^q\colon\U(G_q)\to\U(G)\}_{q\in(0,1]}$ be a continuous family of $*$-isomorphisms with $\varphi^1=\iota$. By definition of $\Gamma_{alg}((\C[G_q])_{q\in[\eps,1]})$ and by linearity we may assume that $a(q)=\hat\varphi^q(a')$ and $b(q)=\hat\varphi^q(b')$ for some $a',b'\in\C[G]$. Let $w_0=s_{i_1}\dots s_{i_n}$ be the longest element in the Weyl group written in reduced form. Consider the homomorphism
$$
\Theta^q\colon C(G_q)\to C(SU_{q^{d_{i_1}}}(2))\otimes\dots\otimes C(SU_{q^{d_{i_n}}}(2)), \ \ \Theta^q(x)=(\sigma^q_{i_1}\otimes\dots\otimes \sigma^q_{i_n})\Delta_q^{(n-1)}(x),
$$
where $\sigma^q_i\colon C(G_q)\to C(SU_{q_i}(2))$ is the $*$-homomorphism which is dual to the embedding $U_{q_i}\sltwo\hookrightarrow U_q\g$ corresponding to the simple root $\alpha_i$. Since $\hat\varphi^q$ are coalgebra maps, we then have

\medskip
$\displaystyle\Theta^q([a(q),b(q)])=\Theta^q([\hat\varphi^q(a'),\hat\varphi^q(b')])$
$$
=\sum_{k=0}^{n-1}\sigma^q_{i_1}(\hat\varphi^q(b'_{(0)})\hat\varphi^q(a'_{(0)}))\otimes\dots\otimes \sigma^q_{i_{k+1}}([\hat\varphi^q(a'_{(k)}),\hat\varphi^q(b'_{(k)})])\otimes\dots\otimes \sigma^q_{i_n}(\hat\varphi^q(a'_{(n-1)})\hat\varphi^q(b'_{(n-1)})),
$$
where we use Sweedler's sumless notation for the coproduct  $\Delta$ on $\C[G]$.
Since $\Delta^{(n-1)}\colon\C[G]\to\C[G]^{\otimes n}$ is a Poisson map with respect to the product Poisson structure on~$\C[G]^{\otimes n}$, we also have
$$
\Theta^q\hat\varphi^q(\{a',b'\})
=\sum_{k=0}^{n-1}\sigma^q_{i_1}\hat\varphi^q(a'_{(0)}b'_{(0)})\otimes\dots\otimes \sigma^q_{i_{k+1}}\hat\varphi^q(\{a'_{(k)},b'_{(k)}\})\otimes\dots\otimes \sigma^q_{i_n}\hat\varphi^q(a'_{(n-1)}b'_{(n-1)}).
$$
By the classification of the irreducible representations of $C(G_q)$ we know that the homomorphism~$\Theta^q$ is an isometry. We thus see that it suffices to show that for any $a,b,c\in\Gamma_{alg}((\C[G_q])_{q\in[\eps,1]})$ with $\{a(1),b(1)\}=c(1)$ and any $1\le j\le r$ we have
$$
\lim_{h\downarrow0}\left\|\sigma_j^{e^{-h}}\left([a(e^{-h}),b(e^{-h})]/{ih}-c(e^{-h})\right)\right\|=0.
$$
Since the family of homomorphisms $(\sigma_j^q)_q$ maps $\Gamma_{alg}((\C[G_q])_{q\in[\eps,1]})$ into $\Gamma_{alg}((\C[SU_{q^{d_j}}(2)])_{q\in[\eps,1]})$ (see the proof of \cite[Proposition 1.2]{NT6}), and $\sigma^1_j\colon \C[G]\to \C[SU(2)]$ is a homomorphism of Poisson algebras, when $SU(2)$ is given the Poisson structure defined by the classical $r$-matrix $i d_j(F\otimes E-E\otimes F)$, to prove the theorem it is therefore enough to consider $G=SU(2)$ with the standard normalization of the invariant form on $\sltwo(\C)$ and the classical $r$-matrix $i(F\otimes E-E\otimes F)$.

The space $\Gamma_{alg}((\C[SU_{q}(2)])_{q\in[\eps,1]})$ is generated as an involutive $C[\eps,1]$-algebra by the sections $q\mapsto\alpha_q$ and $q\mapsto\gamma_q$ (see again the proof of \cite[Proposition 1.2]{NT6}). It follows that it suffices to consider the following four pairs of $(a(q),b(q))$:  $(\alpha_q,\alpha^*_q)$, $(\alpha_q,\gamma_q)$, $(\alpha_q,\gamma^*_q)$ and $(\gamma_q,\gamma^*_q)$. By \eqref{epoissonbracket} the Poisson bracket of $a',b'\in\C[G]$ is given by
$$
\{a',b'\}=(a'_{(1)}\otimes b'_{(1)})(r)a'_{(0)}b'_{(0)}-(a'_{(0)}\otimes b'_{(0)})(r)a'_{(1)}b'_{(1)},
$$
from which we compute
$$
\{\alpha_1,\alpha_1^*\}=-2i\gamma_1\gamma_1^*,\ \
\{\alpha_1,\gamma_1\}=i\alpha_1\gamma_1,\ \
\{\alpha_1,\gamma_1^*\}=i\alpha_1\gamma_1^*,\ \
\{\gamma_1,\gamma_1^*\}=0.
$$
By the relations in $\C[SU_q(2)]$ this gives the result: for instance, $$\frac{[\alpha_{e^{-h}},\alpha_{e^{-h}}^*]}{ih}=\frac{1-e^{-2h}}{ih}\gamma^*_{e^{-h}}\gamma_{e^{-h}}\to -2i\gamma_1^*\gamma_1\ \ \hbox{as}\ \ h\to0.$$
\ep

Recall \cite{Ri}, \cite{La} that a strict quantization of a commutative Poisson $*$-algebra $\A$ is a continuous field $(A_h)_{h\in[0,\delta]}$ of C$^*$-algebras together with a linear map $\QQ=(\QQ_h)_h\colon \A\to\Gamma((A_h)_{h\in[0,\delta]})$ such that~$\A$ is a dense $*$-subalgebra of $A_0$, $\QQ_0=\iota$, $\QQ_h(\A)$ is a dense subspace of $A_h$ for every $h$, and
$$
\lim_{h\downarrow0}\left\|{[\QQ_h(a),\QQ_h(b)]}/{ih}-\QQ_h(\{a,b\})\right\|=0\ \ \hbox{for all}\ \ a,b\in\A.
$$
The pair $((A_h)_h,(\QQ_h)_h)$ is called a strict deformation quantization of $\A$ if in addition every map $\QQ_h$ is injective and its image is a $*$-subalgebra of $A_h$.

The structure which emerges from Theorem~\ref{tquantization} is only slightly different: we have a continuous field $(A_h)_{h\in[0,\delta]}$ of C$^*$-algebras and a dense involutive $C[0,\delta]$-subalgebra $\QQ$ of $\Gamma((A_h)_{h\in[0,\delta]})$ such that~$\A$ is a dense $*$-subalgebra of $A_0$, the image of $\QQ$ in $A_0$ coincides with $\A$ and
$$
\lim_{h\downarrow0}\frac{[a(h),b(h)]}{ih}=\{a(0),b(0)\}\ \ \hbox{for all}\ \ a,b\in\QQ.
$$
The advantage of this formulation is that in our examples this structure is completely canonical. If one however insists on the standard formulation of deformation quantization, the required maps $\QQ_h$ are in abundance, but it is impossible to make a canonical choice.

\begin{corollary}
Let  $\{\varphi^q\colon\U(G_q)\to\U(G)\}_{q\in(0,1]}$ be a continuous family of $*$-isomorphisms with $\varphi^1=\iota$ such that $\varphi^q(\U(K^{S,L}_q))=\U(K^{S,L})$. Then the pair $((A_h)_{h\in[0,\delta]},(\QQ_h)_{h\in[0,\delta]})$, where $A_h=C(G_{e^{-h}}/K^{S,L}_{e^{-h}})$ and $\QQ_h$ is the restriction of $\hat\varphi^{e^{-h}}$ to $\C[G/K^{S,L}]$, defines a strict deformation quantization of the Poisson algebra $\C[G/K^{S,L}]$ for any $\delta>0$.
\end{corollary}

\begin{remark}
Sometimes one also requires the maps $\QQ_h$ to be $*$-preserving. The maps in the above corollary do not satisfy this property, but it is easy to modify them to get maps that are $*$-preserving. To do this, for every $\lambda\in P_+$ and $h\ge0$ consider the subspace $\tilde\A_h^\lambda\subset \C[G_{e^{-h}}]$ spanned by the matrix coefficients of the irreducible representations with highest weights $\lambda$ and $-w_0\lambda$. Put $\A_h^\lambda= \tilde \A_h^\lambda\cap \C[G_{e^{-h}}/K^{S,L}_{e^{-h}}]$. Then $\A_h^\lambda$ is a finite dimensional selfadjoint subspace of $A_h$ and $\QQ_h$ maps $\A_0^\lambda$ onto $\A_h^\lambda$. Then $R^\lambda_h=\QQ^{-1}_h((\A^\lambda_h)_{sa})\subset \A_0^\lambda$ is a continuous family of real forms of the space $\A^\lambda_0$. Hence there exists a continuous family of linear isomorphisms $T^\lambda_h\colon \A_0^\lambda\to \A^\lambda_0$ such that $T^\lambda_h(R^\lambda_0)=R^\lambda_h$ and $T^\lambda_0=\iota$. The space $\C[G/K^{S,L}]$ is the direct sum of the spaces $\A_0^\lambda$ over a set of representatives $\lambda$ of the quotient space of $P_+$ by the action of the involution $-w_0$. Fixing such a direct sum decomposition define $T_h\colon \C[G/K^{S,L}]\to\C[G/K^{S,L}]$ using the operators $T_h^\lambda$. Then the maps $\QQ_hT_h\colon\C[G/K^{S,L}]\to \C[G_{e^{-h}}/K^{S,L}_{e^{-h}}]$ are $*$-preserving linear isomorphisms defining a strict deformation quantization of $\C[G/K^{S,L}]$.
\end{remark}

\bigskip


\section{K-theory} \label{s6}

In this section we will show that the C$^*$-algebras $C(G_q/K^{S,L}_q)$ are KK-equivalent to $C(G/K^{S,L})$. In fact we will prove the following more precise and stronger result, which is important for applications~\cite{NT6}.

\begin{theorem} \label{tNagy}
For any $0<a<b\le1$ and $q_0\in[a,b]$ the evaluation map $$\Gamma((C(G_q/K^{S,L}_q))_{q\in[a,b]})\to C(G_{q_0}/K^{S,L}_{q_0})$$ is a KK-equivalence.
\end{theorem}

That $C(SU_q(N))$ is KK-equivalent to $C(SU(N))$, was proved by Nagy~\cite{Na} using the composition series obtained by Sheu~\cite{Sh}. In view of Theorem~\ref{tcompseries} the general case of $C(G_q/K^{S,L}_q)$ is virtually the same, but since it might seem that the proof of Nagy depends in an essential way on the extension of E-theory developed in~\cite{Na2}, we will give a complete argument within just the standard KK-theoretic framework.

\smallskip

We will repeatedly use the following basic properties of KK-equivalence.

\begin{lemma} Let $0\to J\to A\xrightarrow{\pi} A/J\to0$ be a semisplit short exact sequence of separable C$^*$-algebras. Then the following conditions are equivalent:
\enu{i} the homomorphism $\pi\colon A\to A/J$ is a KK-equivalence;
\enu{ii} the map $\pi_*\colon KK(D,A)\to KK(D,A/J)$ is an isomorphism for every separable C$^*$-algebra $D$;
\enu{iii} the map $\pi^*\colon KK(A/J,D)\to KK(A,D)$ is an isomorphism for every separable C$^*$-algebra $D$;
\enu{iv} the C$^*$-algebra $J$ is KK-contractible, that is, $KK(J,J)=0$.
\end{lemma}

\bp Equivalence of (ii), (iii) and (iv) follows from the two $6$-term exact sequences in KK-theory associated with $0\to J\to A\xrightarrow{\pi} A/J\to0$. That (i) implies (ii) is immediate. Finally, that (ii) implies~(i) follows from the general observation that if $f\colon X\to Y$ is a morphism in some category $\CC$ such that for every object $Z$ the map $\operatorname{Mor}(Z,X)\to\operatorname{Mor}(Z,Y)$, $g\mapsto f\circ g$, is a bijection, then $f$ is an isomorphism.
\ep

Note that all the algebras appearing in this section will be of type I, hence nuclear, so all the short exact sequences will automatically be semisplit.

\smallskip

To prove the theorem we will first establish the analogous result for quantum disks.

\begin{lemma} \label{ldisk}
For any $0\le a<b\le1$ and $q_0\in[a,b]$ the evaluation maps $\Gamma((C_0(\DD_q))_{q\in [a,b]})\to C_0(\DD_{q_0})$  and $\Gamma((C(\bar\DD_q))_{q\in [a,b]})\to C(\bar\DD_{q_0})$ are KK-equivalences.
\end{lemma}

\bp Any of the two $6$-term exact sequences in KK-theory applied to the exact rows of the commutative diagram
$$
\xymatrix{
0\ar[r] & \Gamma((C_0(\DD_q))_{q\in [a,b]})\ar[d]\ar[r] & \Gamma((C(\bar\DD_q))_{q\in [a,b]})\ar[r]\ar[d] & C[a,b]\otimes C(\T)\ar[r]\ar[d] & 0\\
0\ar[r] & C_0(\DD_{q_0}) \ar[r] & C(\bar\DD_{q_0})\ar[r] & C(\T)\ar[r] & 0
}
$$
implies that it suffices to show that $\Gamma((C(\bar\DD_q))_{q\in [a,b]})\to C(\bar\DD_{q_0})$ is a KK-equivalence. Observe also that for $q_0\in(a,b)$ the kernel of $\Gamma((C(\bar\DD_q))_{q\in [a,b]})\to C(\bar\DD_{q_0})$ is the direct sum of the kernels of $\Gamma((C(\bar\DD_q))_{q\in [a,q_0]})\to C(\bar\DD_{q_0})$ and $\Gamma((C(\bar\DD_q))_{q\in [q_0,b]})\to C(\bar\DD_{q_0})$. So to prove that the kernels of the evaluation maps are KK-contractible it suffices to consider the evaluations at the end points.

Since the field $(C(\bar\DD_q))_{q\in[0,1)}$ is constant with fiber $\TT$, the kernel of $\Gamma((C(\bar\DD_q))_{q\in [a,b]})\to C(\bar\DD_b)$ is isomorphic to $C_0[a,b)\otimes\TT$, hence it is contractible. Similarly, if $b<1$ then the kernel of $\Gamma((C(\bar\DD_q))_{q\in [a,b]})\to C(\bar\DD_{a})$ is contractible.

It remains to prove that $ev_a\colon \Gamma((C(\bar\DD_q))_{q\in [a,1]})\to C(\bar\DD_{a})$ is a KK-equivalence. Since we already know that $\Gamma((C(\bar\DD_q))_{q\in [a,1]})\to C(\bar\DD_{1})=C(\bar\DD)$ is a KK-equivalence, the C$^*$-algebra $\Gamma((C(\bar\DD_q))_{q\in [a,1]})$ is KK-equivalent to~$\C$ and the group $K_0(\Gamma((C(\bar\DD_q))_{q\in [a,1]}))$ is generated by $[1]$. But it is also well-known that the C$^*$-algebra $C(\bar\DD_a)\cong\TT$ is KK-equivalent to $\C$ and its $K_0$-group is generated by~$[1]$. Therefore we just have to check that the KK-class of $ev_a$ is a generator of
$$
KK(\Gamma((C(\bar\DD_q))_{q\in [a,1]}),C(\bar\DD_a))\cong KK(\C,\C)\cong\Z.
$$
Since ${ev_a}_*\colon K_0(\Gamma((C(\bar\DD_q))_{q\in [a,1]}))\to K_0(C(\bar\DD_a))$ is an isomorphism, this is clearly the case.
\ep

We remark that the last part of the above proof can be slightly shortened by using the Universal Coefficient Theorem. Similarly the next lemma can be quickly deduced from Lemma~\ref{ldisk} using Kasparov's ${\mathcal R}$KK-groups. Since both proofs are quite short anyway, we prefer to keep things as elementary as possible.

\begin{lemma} \label{lpolydisk}
Assume $p_1,\dots,p_n>0$ and $0\le a<b\le1$.
Then for any $q_0\in[a,b]$ the evaluation map
$$
\Gamma((C_0(\DD_{q^{p_1}})\otimes\dots\otimes C_0(\DD_{q^{p_n}}))_{q\in[a,b]})\to C_0(\DD_{q^{p_1}_0})\otimes\dots\otimes C_0(\DD_{q^{p_n}_0})
$$
is a KK-equivalence.
\end{lemma}

Here the family $(C_0(\DD_{q^{p_1}})\otimes\dots\otimes C_0(\DD_{q^{p_n}}))_{q\in[a,b]}$ is of course given the unique continuous field structure such that the tensor product of continuous sections is a continuous section. That such a structure exists, can be checked by the same argument as in the proof of Lemma~\ref{ldiskcont}, but this is also a consequence of a general result of Kirchberg and Wassermann \cite[Theorem~4.6]{KW} saying that if $(A_q)_q$ and $(B_q)_q$ are continuous fields of C$^*$-algebras and $\Gamma((A_q)_q)$ is exact then $(A_q\otimes B_q)_q$ is a continuous field.

\bp[Proof of Lemma~\ref{lpolydisk}] To simplify the notation assume $p_1=\dots=p_n=1$.
The proof of the lemma is by induction on $n$. Furthermore, it is convenient to simultaneously prove the same result for the continuous fields of the C$^*$-algebras
$$
A_{m,n}^q=\underbrace{C_0(\DD_q)\otimes\dots\otimes C_0(\DD_q)}_m\otimes \underbrace{C(\bar\DD_q)\otimes\dots\otimes C(\bar\DD_q)}_{n-m}
$$
for $m=0,\dots,n$. For $n=1$ the result is proved in Lemma~\ref{ldisk}. Assume $n>1$. We will prove that the evaluation map $\Gamma((A^q_{m,n})_{q\in[a,b]})\to A^{q_0}_{m,n}$ is a KK-equivalence by induction on $m$. For $m=0$ the proof is literally the same as that of Lemma~\ref{ldisk}, with $\TT$ replaced by $\TT^{\otimes n}$. For $m\ge1$ applying $\otimes A^q_{m-1,n-1}$ to the exact sequence $0\to C_0(\DD_q)\to C(\bar\DD_q)\to C(\T)\to 0$ we get an exact sequence
\begin{equation} \label{eshortdisks}
0\to A^q_{m,n}\to A^q_{m-1,n}\to A^q_{m-1,n-1}\otimes C(\T)\to 0.
\end{equation}
Since the evaluation map $\Gamma((A^q_{m-1,n-1})_{q\in[a,b]})\to A^{q_0}_{m-1,n-1}$ is a KK-equivalence by the inductive assumption on $n$, the map
$$
\Gamma((A^q_{m-1,n-1}\otimes C(\T))_{q\in[a,b]})=\Gamma((A^q_{m-1,n-1})_{q\in[a,b]})\otimes C(\T)\to A^{q_0}_{m-1,n-1}\otimes C(\T)
$$
is a KK-equivalence as well. Since $\Gamma((A^q_{m-1,n})_{q\in[a,b]})\to A^{q_0}_{m-1,n}$ is a KK-equivalence by the inductive assumption on $m$, applying one of the $6$-term exact sequences in KK-theory to \eqref{eshortdisks} we conclude that $\Gamma((A^q_{m,n})_{q\in[a,b]})\to A^{q_0}_{m,n}$ is also a KK-equivalence.
\ep

\bp[Proof of Theorem~\ref{tNagy}]
Consider the ideals $J_m^q\subset C(G_q/K^{S,L}_q)$ defined in Theorem~\ref{tcompseries}. Since they are the fiber-wise kernels of morphisms of continuous fields of C$^*$-algebras, they form continuous subfields of C$^*$-algebras of $(C(G_q/K^{S,L}_q))_q$, see e.g.~\cite[Proposition~2.6(ii)]{Na}.  Furthermore, we have short exact sequences $0\to J^q_{m}\to J^q_{m-1}\to A^q_m\to0$ and corresponding short exact sequences of the C$^*$-algebras of continuous sections, where
$$
A_m^q=\bigoplus_{w\in W^S: \ell(w)=m}C(T/T_L)\otimes C_0(\DD_{q^d{_{i_1(w)}}})\otimes\dots\otimes
C_0(\DD_{q^{d_{i_m(w)}}}).
$$
By Lemma~\ref{lpolydisk} the evaluation maps $\Gamma((A^q_m)_q)\to A^m_{q_0}$ are KK-equivalences. As $J^q_{m_0}=0$, using the $6$-term exact sequences in KK-theory we prove that $\Gamma((J^q_m)_q)\to J_m^{q_0}$ are KK-equivalences for all~$m$ by downward induction from $m=m_0$ to $m=-1$. The case $m=-1$ is the statement of the theorem.
\ep

Since the continuous field structure on $(C(G_q/K^{S,L}_q))_q$ does not depend on any choices, we therefore know that the K-groups of $C(G_q/K^{S,L}_q)$ are canonically isomorphic to those of $C(G/K^{S,L})$, but this gives no information about explicit generators of these groups. Some information can however be extracted. Let $e$ be a projection in a matrix algebra over $C(G/K^{S,L})$. Assume we can find a continuous field of projections $e(q)$ in matrix algebras over $C(G_q/K^{S,L}_q)$ such that $e(1)=e$. Then the class of $e(q)$ in $K_0(C(G_q/K^{S,L}_q))$ is exactly the class corresponding to $[e]$ under the KK-equivalence between $C(G_q/K^{S,L}_q)$ and $C(G/K^{S,L})$.

As a simple example consider the Podle\'s sphere $S^2_q=SU_q(2)/\T$. It is well-known, and follows immediately from Theorem~\ref{tcompseries}, that the homomorphism $\rho_q$ defines an isomorphism of $C(S^2_q)$ onto the unitization $C_0(\DD_q)^\sim\subset C(\bar\DD_q)$ of $C_0(\DD_q)$. So for $q\in(0,1)$ the C$^*$-algebra $C(S^2_q)$ is isomorphic to the algebra of compact operators on $\ell^2(\Z_+)$ with unit adjoined. From this point of view the most natural generators of $K_0(C(S^2_q))\cong\Z^2$ are $[1]$ and the class of the rank-one projection onto $\C e_0$. The latter projection has no meaning for $q=1$. On the other hand, $K_0(S^2)$ is generated by $[1]$ and the class of the Bott element. Under the identification $S^2=SU(2)/\T$ this class can be represented by the projection $\begin{pmatrix}\gamma^*_1\gamma_1 & -\alpha_1\gamma_1^*\\ -\gamma_1\alpha_1^* & \alpha_1^*\alpha_1\end{pmatrix}$. This projection belongs to the continuous family of projections $e(q)=\begin{pmatrix}q^2\gamma^*_q\gamma_q & -\alpha_q\gamma_q^*\\ -\gamma_q\alpha_q^* & \alpha_q^*\alpha_q\end{pmatrix}$, see~\cite{HM}. Therefore $K_0(C(S^2_q))$ is generated by $[1]$ and $[e(q)]$.

As another example, from the classical result of Hodgkin~\cite{Ho} we conclude that the fundamental corepresentations of $C(G_q)$ define independent generators of $K_1(C(G_q))$ (but not all of them if the rank of $G$ is at least $3$).

It would be interesting to develop a general technique for how to lift K-theory classes for $G/K^{S,L}$ to $\Gamma((C(G_q/K^{S,L}_q))_q)$.

\end{document}